\documentclass[a4 paper,12pt]{article}
\usepackage[T1]{fontenc}
\usepackage[latin1]{inputenc}
\usepackage[german,english]{babel}
\usepackage{graphicx}
\usepackage{mathrsfs}
\usepackage{amsfonts}
\usepackage{pifont}
\usepackage{amsxtra}
\usepackage{amssymb}
\usepackage{eufrak}
\usepackage{pst-all}
\usepackage{multido}
\usepackage{amsmath}
\usepackage{color}
\usepackage{makeidx}
\usepackage{multicol}
\usepackage[absolut]{overpic}

\usepackage[rm,bf,tiny,center]{titlesec}
\titlelabel{\thetitle.\enspace}
\newcommand{\ger}[1]{\mathfrak{#1}}

\newcommand{\OO}{\mathcal{O}}

\newcommand{\sur}{\twoheadrightarrow}
\newcommand{\A}{\mathcal{A}}
\newcommand{\ma}{\leqslant}
\newcommand{\De}{\Delta}
\newcommand{\CC}{\mathbb{C}}

\newcommand{\DD}{\mathcal{D}}
\newcommand{\II}{\mathcal{I}}

\DeclareMathOperator{\SM}{SM}
\DeclareMathOperator{\add}{add}
\DeclareMathOperator{\dimm}{dim}
\DeclareMathOperator{\Hom}{Hom}
\DeclareMathOperator{\Kern}{ker}
\DeclareMathOperator{\Loc}{Loc}
\DeclareMathOperator{\End}{End}
\DeclareMathOperator{\topp}{top}
\DeclareMathOperator{\rad}{rad}
\DeclareMathOperator{\modd}{mod}
\DeclareMathOperator{\soc}{soc}

\DeclareMathOperator{\Ext}{Ext}
\DeclareMathOperator{\im}{im}
\DeclareMathOperator{\spann}{span}

\newtheorem{num}{}[section]

\usepackage{geometry}

\usepackage{fancyhdr}

\begin{document}




\begin{center}
\begin{large}\textbf{1-quasi-hereditary algebras} \end{large}
\end{center}
\begin{center}
Daiva Pu\v{c}inskait\.{e}
\end{center}

\begin{abstract}
Motivated by the structure of the  algebras associated to the blocks of the BGG-category $\OO$, we define a subclass of quasi-hereditary algebras called 1-quasi-hereditary. 
Many properties of these  algebras only depend  on the defining  partial order.   In particular,  we can determine the quiver and  the form of the relations.   Moreover, if the Ringel dual  of a  1-quasi-hereditary algebra   is also 1-quasi-hereditary, then the structure of the characteristic tilting module can be computed. 
\end{abstract}

\begin{center}
\textbf{Introduction}
\end{center}
  
The class of quasi-hereditary algebras,  
 defined  by   Cline, Parshall and Scott  \cite{CPS}, can be   regarded as a 	generalization  of the algebras  associated  to  the blocks  of the  Bernstein-Gelfand-Gelfand  category $\OO(\ger{g})$  of  a complex semisimple  Lie algebra  $\ger{g}$  (see \cite{BGG}).  Every block $\mathcal{B}(\ger{g})$ is equivalent to the  category of modules
over a finite dimensional   $\CC$-algebra  $\A_{\mathcal{B}}(\ger{g})$.

The algebras  $\A_{\mathcal{B}}(\ger{g})$  are   BGG-algebras as defined in \cite{Ir} and in  \cite{Chan}. They are endowed    with  a duality functor  on their module category  which fixes the simple modules. Another important structural feature is the presence  of exact Borel subalgebras and $\De$-subalgebras introduced by König in \cite{K}. These subalgebras 
 provide  a correspondence between   $\De$-good filtrations  and  Jordan-Hölder-filtrations. Moreover, Soergel has shown that  $\A_{\mathcal{B}}(\ger{g})$ is Morita equivalent to its Ringel dual $R(\A_{\mathcal{B}}(\ger{g}))$ (see \cite{So}).
 

Motivated by these results,   
      in this paper we introduce  a   class of    quasi-hereditary algebras,  called 1-quasi-hereditary.  Among other  properties they are characterized by the fact that  
      all possible non-zero filtration-multiplicities for  $\De$-good filtrations of indecomposable projectives  and Jordan-Hölder filtrations of standard modules  are equal to 1.

  The class of 1-quasi-hereditary algebras  is   related to the  aforementioned   classes of  quasi-hereditary algebras:  
   Many   factor algebras (related to  saturated subsets)  of an  algebra  of type  $\A_{\mathcal{B}}(\ger{g})$ are 1-quasi-hereditary.   The understanding of 1-quasi-hereditary algebras gives some information on  the relations,  the structure of the characteristic tilting module etc. of  $\A_{\mathcal{B}}(\ger{g})$.
 Another class of examples is provided by the  quasi-hereditary algebras considered by    Dlab, Heath and Marko in \cite{DHM}. These  algebras   are 1-quasi-hereditary BGG-algebras, however  1-quasi-hereditary algebras are in general not BGG-algebras.   All known 1-quasi-hereditary algebras have  exact Borel  and $\De$-subalgebras.  Several  examples, which  show  the complexity of such   algebras  and  their  additional properties are presented in \cite{P}.

Our first main result  shows   that many invariants of 1-quasi-hereditary algebras  depend  only on the given  partial ordering:
\\

\hspace*{-5mm}\textbf{Theorem A.} \textit{Let $A=(KQ/\II,\leqslant)$  be a (basic)  1-quasi-hereditary algebra.
  Then  
\begin{itemize}
	\item[(1)] $Q$  is the  double    of the  quiver of the  incidence algebra corresponding to  $\ma$, \\[3pt]
	i.e. $Q_1=\left\{\right.$\psset{xunit=1mm,yunit=1mm,runit=1mm}
\begin{pspicture}(0,1)(0,0)
\begin{small}
\rput(2.5,1){
\rput(0,0){\rnode{0z}{$i$}}
\rput(10,0){\rnode{1z}{$j$}}
}
\psset{nodesep=2pt,offset=2pt,arrows=<-}
\ncline{0z}{1z}
\ncline{1z}{0z} \end{small}
\end{pspicture} \hspace{1.5cm}  $\mid  \ \  i$ and $j$ are neighbours w.r.t. $\ma$ $\left.\right\}$.  
\\[30pt]
 \hspace*{-0.8cm}\begin{minipage}{6cm}
 \hrule
 \vspace{0.5cm} \textcolor{black}{ }
\end{minipage}
\\[-4mm]
\hspace*{-0.7cm}\begin{footnotesize}\rm{Partly supported by the D.F.G. priority program SPP 1388 ``Darstellungstheorie''.}\end{footnotesize}
	\item[(2)]  $\II$ is generated by the relations of the form  $\displaystyle p-\!\sum_{j,k\ma i} c_i\cdot p(j,i,k)  $, where $p=(j\!\to\!\cdots\!\to\!k)$   and   $p(j,i,k)$  are paths in $Q$  of  the   form  $ (j\!=\!j_1\! \to\!\cdots\! \to \!j_m\!\to\! i \!\to\! k_1 \!\to \!\cdots\!\to \!k_r\!=\!k)$ with $j_1\!<\!\cdots\! < \!j_m\!<i>\!k_1\!>\!\cdots \!>k_r$.
		\item[(3)]  The   $\De$-good filtrations  of the projective indecomposable module at the vertex $i\in Q_0$ are in one-to-one correspondence with special sequences of   vertices $j$ with $ j \geqslant i$. 
\end{itemize}}

  An important feature in  the representation theory of quasi-hereditary algebras is the concept of the Ringel dual: The algebra $R(A):=\End_A(T)^{op}$ is quasi-hereditary, where $T=\bigoplus_{i\in Q_0}T(i)$ is the    characteristic tilting module. 
In view of Soergels work, this raises the question whether  the class of 1-quasi-hereditary algebras is  closed under Ringel-duality. 
\\

\hspace*{-5mm}\textbf{Theorem B}. \textit{Let $A=(KQ/\II, \ma)$ be a   1-quasi-hereditary algebra.   Then 
\begin{center}
$R(A)$  is 1-quasi-hereditary \  if and only if \  $T(i)$ is local for any $i\in Q_0$. 
\end{center}}
Moreover,  in this case  we   have  a precise description of  $T(i)$.
\\

Our paper is 	organised  as follows: In Section  1, we introduce some notations, recall some definitions and  basic facts  for later use.

 In Section   2,  we give several   properties  of 1-quasi-hereditary algebras, which can be derived from the definition using the general representation  theory of bound quiver algebras.  These  properties are essential for the proof of Theorem A (1).  
 
In Section  3, we  present a  particular   
 basis of a 1-quasi-hereditary algebra $A$, which  can be described combinatorially and only  depends  on the corresponding partial order (it  consist the  paths of the form $p(j,i,k)$).   Consequently,  we obtain  a 	system   of  relations of $A$ described in Theorem A (2).

In   Section 4,  we  determine the  set of $\De$-good  filtrations of all projective  indecomposable modules over 1-quasi-hereditary algebras and establish their  	relationship with the Jordan-Hölder-filtrations of costandard modules. Using the result of Ringel \cite{Rin1}, which says  that the subcategory $\ger{F}(\De)$  is resolving,   we  determine all local  modules having $\De$-good  filtrations. 
We also record the dual results.  

In  Section 5,  we consider    factor algebras $A(i):=A/A(\sum_{j\not\ma i}e_j)A$  for   $i\in Q_0$ of a 1-quasi-hereditary algebra $A$, where $e_i$ is a  primitive  idempotent. 
If $A(i)$ is 1-quasi-hereditary, then we obtain an explicit   expression of the direct summand $T(i)$ of the characteristic tilting module. 

Using these results  in Section 6,   we turn to the question when the Ringel dual of a 1-quasi-hereditary algebra is 1-quasi-hereditary. We elaborate on Theorem B by establishing necessary and sufficient conditions involving the structure of tilting modules and projective indecomposable modules.

\section{Preliminaries}
\begin{small}
Throughout the paper, $\A$  denotes a finite dimensional, basic  $K$-algebra over an algebraically closed field $K$,  which will be   represented  by a quiver and relations (Theorem of Gabriel) and $\modd \A$ is the category of finite dimensional left $\A$-modules. 
 In the following part   we will focus on  some general facts from the representations theory of bound quiver algebras, which we will use  in this paper. The relevant material can be found in  \cite{ASS}. 
 \end{small}
\\

We consider   algebras   $\A=KQ/\II$  and by $Q_0$ (resp. $Q_1$) we denote  the set of  vertices  (resp. the set of arrows) in $Q$. 
For any  $i\in Q_0$ the corresponding  trivial path will be denoted    by $e_i$,  the   simple module, the projective indecomposable  and the  injective  indecomposable  $\A$-module, will be  denoted by  $S(i)$,  $P(i)$ and  $I(i)$ respectively. 
A  path  $p=\left(j\rightarrow \cdots \rightarrow i \rightarrow \cdots \rightarrow k\right)$  
 is the product   of paths 
   $p_1=\left(i \rightarrow \cdots \rightarrow k\right)$ and $p_2=\left(j\rightarrow \cdots \rightarrow i\right)$   written as  $p=p_1\cdot p_2$. The   $\A$-map   corresponding to  $p $  is given by  $f_{p}: P(k)\rightarrow P(j)$ via  $f_{p}(a\cdot e_k)= a\cdot p\cdot e_j$  for all  $a\in \A$ and we have  $f_{p}= f_{p_2}\circ f_{p_1}$.

 For any   $M\in \modd \A$  it is  $M\cong \bigoplus_{i\in Q_0}M_i$, where $M_i$ is the  subspace of $M$ corresponding to  $i\in Q_0$.   
  We denote by $\left[M:S(i)\right]=\dimm_K M_i$  the Jordan-Hölder multiplicity of $S(i)$ in $M$.
 For any $m\in M$,  we denote   by  $\left\langle m\right)$ the submodule of $M$ generated by $m$  (i.e. $\left\langle m\right)=\A \cdot m$). The set of all local submodules  of $M$ with  $\topp$ isomorphic to $S(i)$,  we  denote  by  $\Loc_i(M)$.  It is clear that        $\Loc_i(M) =\left\{\left\langle m\right)\mid m\in M_i \backslash \left\{0\right\}\right\}=  \left\{\im (f) \mid f\in \Hom_A(P(i),M), \ f\neq 0\right\}$. 
\\

The definition of quasi-hereditary  algebras introduced by Cline-Parshall-Scott \cite{CPS}
  implies in  particular  the presence   of a   partial order   on the vertices of the corresponding  quiver.  The equivalent definition and relevant terminology is given by Dlab and Ringel in \cite{DRin}.  To recap briefly:  For an algebra  $\A\cong KQ/\II$ let  $(Q_0,\ma)$ be a partially ordered set.     
For  every   $i\in Q_0$ the  \textit{standard}  module    $\De(i)$  is   the largest  factor module    of 
    $P(i)$  such that      $[\De(i):S(j)]=0$ for all $j\in Q_0$ with  $j\not\ma i$  and    the  \textit{costandard}  module   $\nabla(i)$ is  the largest submodule of 
    $I(i)$ such that   $[\nabla(i):S(j)]=0$ for all $j\in Q_0$ with   $j\not\ma i$.   We denote   by $\ger{F}(\De)$
       the full subcategory of $\modd \A$  consisting of the modules having a filtration  such that each subquotient is isomorphic to a standard module.    The modules   in $\ger{F}(\De)$  are called    \textit{$\De$-good}  and  the corresponding    filtrations are   \textit{$\De$-good filtrations}      (resp.  \textit{$\nabla$-good} modules have     \textit{$\nabla$-good filtrations} and belong  to $\ger{F}(\nabla)$).   For  $M\in \ger{F}(\De)$,   we   denote  by $\left(M:\De(i)\right)$  the  (well-defined)  number of  subquotients  isomorphic to  $\De(i)$ in some $\De$-good filtration  of $M$ (resp. $\nabla(i)$ appears  $(M:\nabla(i))$ times   in some  $\nabla$-good filtration of $M\in \ger{F}(\nabla)$). 
\\[8pt]
\hspace*{2mm}\parbox{16cm}{
The algebra $\A=\left(KQ/\II, \ma\right)$ is   \textit{quasi-hereditary}    if  for all $i,j\in Q_0$ the following   holds: 
\begin{itemize}
	\item   $[\De(i):S(i)]=1$,
	\item $P(i)$ is a  $\De$-good module   with $\left(P(i):\De(j)\right)=0$ for all $j\not\geqslant i$ and $\left(P(i):\De(i)\right)=1$.
\end{itemize}
     }

\begin{num}\normalfont{\textbf{Remark.}}  If  $(\A,\ma)$ is quasi-hereditary, then for any $i\in Q_0(\A)$ the following  holds: 
\\[5pt]
\hspace*{0.7cm}$\displaystyle \Delta(i)= P(i)/\left(\sum_{i < j}\sum_{f\in \Hom_A(P(j),P(i))}\im (f)\right)$    \ \  resp.  \ \  $\displaystyle \nabla(i)= \bigcap_{i<
j}\bigcap_{f\in \Hom_A(I(i),I(j))} \Kern (f)$ 
\\[5pt]
Moreover, if $i\in Q_0$ is minimal 	with respect to $\ma$, then   $\De(i)\cong \nabla(i)\cong S(i)$  and if $i\in Q_0$ is maximal  then  $P(i)\cong \De(i)$ as well as $I(i)\cong \nabla(i)$.
\label{qhstandq} \end{num}

\begin{num}\normalfont{\textbf{Definition.}}
A quasi-hereditary algebra $A=(KQ/\II, \leqslant)$
is called 
\textit{1-quasi-hereditary}
if for all $i,j\in Q_0=\left\{1, \ldots , n\right\}$ the
following conditions are satisfied:
\begin{itemize}
    \item[(1)] There is a smallest and a largest element with respect to  $\ma$, 
    \\ without loss of generality we will assume  them to be $1$ resp. $n$,
    \item[(2)] $[\Delta(i):S(j)]=\big(P(j):\Delta(i)\big)=1$ for $j\leqslant i$,  
\item[(3)] $\soc P(j) \cong \topp I(j)\cong  S(1)$,
\item[(4)] $\Delta(i) \hookrightarrow \Delta(n)$ and $\nabla(n)\twoheadrightarrow
\nabla(i). $
\end{itemize} \label{def1qh}
\end{num}
\text{  }\\[-15mm]

The class of 1-quasi-hereditary algebras are related to  several  subclasses of quasi-hereditary algebras:  
 Many   factor algebras (related to a  saturated subsets)  of an algebra  associated  to a block  of the  category $\OO(\ger{g})$  of  a  semisimple  $\CC$-Lie algebra  $\ger{g}$ are 1-quasi-hereditary. 
 If  $\text{rank}(\ger{g})\leq 2$, then  an  algebra corresponding to a block  of $\OO(\ger{g})$ is 1-quasi-hereditary.   
This algebras are   BGG-algebras in the sense of \cite{Chan}  and Ringel self-dual, 	however 1-quasi-hereditary  algebras are not BGG-algebras in general and the class of 1-quasi-hereditary algebras is not closed under Ringel duality.  Moreover all known  1-quasi-hereditary algebras have exact  Borel and  $\De$-subalgebras in sense of König \cite{K}. In \cite{P} we give several  examples to  illustrate this specific properties.
\\
\\
\parbox{14cm}{Let $(Q_0,\ma)$ be the corresponding poset  of a 1-quasi-hereditary algebra $KQ/\II$.  For any $j\in Q_0$,  we define 
 \begin{center}
$\Lambda_{(j)}:= \left\{i\in Q_0\mid i\leqslant j\right\}$ \  \  and \  \   $\Lambda^{(j)}:= \left\{i\in Q_0\mid i\geqslant  j\right\}$
\end{center}
If $i<j$  (resp. $i>j$) and they are neighbours  with respect to $\ma$, then we write 
} \psset{xunit=0.3mm,yunit=0.3mm,runit=1mm}
\begin{pspicture}(50,0)(0,0)
\rput(40,14){
\rput(0,0){\begin{tiny}
\rput(0,15.7){\rnode{n}{\begin{footnotesize}$n$\end{footnotesize}}}
\rput(0,10){\rnode{A}{$\bullet$}}
\rput(0,-60){\rnode{Aa}{$\bullet$}}
\rput(0,-65.9){\rnode{1}{\begin{footnotesize}$1$\end{footnotesize}}}
\rput(0,-60){\rnode{B}{$\bullet$}}
\rput(21,-17){\rnode{aa}{}}
\rput(-21,-17){\rnode{bb}{}}
\rput(-13,2){\rnode{A1}{}}
\rput(-13,-52){\rnode{B1}{}}
\rput(13,2){\rnode{A2}{}}
\rput(13,-52){\rnode{B2}{}}
              \rput(0,-25){\rnode{2}{$\bullet$}}
              
              \rput(8,-10){\rnode{1r}{}}
              \rput(-12,-10){\rnode{2l}{}}
              \rput(10,-25){\rnode{j}{\begin{footnotesize}$j$\end{footnotesize}}}
              \rput(6,-10){\rnode{L}{}}
              \rput(30,0){\rnode{LL}{}}
              \rput(6,-45){\rnode{Lu}{}}
              \rput(30,-35){\rnode{LLu}{}}
            \begin{small} 
             \rput(40,1){\rnode{LLt}{$\Lambda^{(j)}$}}
            \rput(40,-37){\rnode{LLtu}{$\Lambda_{(j)}$}}
            \end{small}
              \end{tiny}
              \ncarc[arcangle=-25,linecolor=black]{->}{LL}{L}
              \ncarc[arcangle=-25,linecolor=black]{->}{LLu}{Lu}
        \ncarc[arcangle=65,linestyle=dashed]{A}{2}
        \ncarc[arcangle=-65,linestyle=dashed]{A}{2}
        \ncarc[arcangle=65,linestyle=dashed]{Aa}{2}
        \ncarc[arcangle=-65,linestyle=dashed]{Aa}{2}
             \psset{nodesep=1pt}
              \ncline[linestyle=dotted]{21}{22}
              \ncline{-}{A}{A1}
              \ncline{-}{A1}{A}
              \ncline{-}{B1}{B}
              \ncline{-}{B}{B1}
              \ncline{-}{A}{A2}
              \ncline{-}{A2}{A}
              \ncline{-}{B2}{B}
              \ncline{-}{B}{B2} 
              \ncarc[arcangle=25,linestyle=dotted,linecolor=black]{A2}{aa}
              \ncarc[arcangle=-25,linestyle=dotted,linecolor=black]{-}{aa}{A2}
              \ncarc[arcangle=25,linestyle=dotted,linecolor=black]{-}{aa}{B2}
              \ncarc[arcangle=-25,linestyle=dotted,linecolor=black]{B2}{aa}
              \ncarc[arcangle=-25,linestyle=dotted,linecolor=black]{A1}{bb}
              \ncarc[arcangle=25,linestyle=dotted,linecolor=black]{-}{bb}{A1}
              \ncarc[arcangle=-25,linestyle=dotted,linecolor=black]{-}{bb}{B1}
              \ncarc[arcangle=25,linestyle=dotted,linecolor=black]{B1}{bb}
\multiput(-1.25,2)(1.5,0){2}{\color{black!70}{\circle*{0.1}}}
\multiput(-2,1)(1.5,0){3}{\color{black!70}{\circle*{0.1}}}
\multiput(-1.25,0)(1.5,0){3}{\color{black!70}{\circle*{0.1}}}
\multiput(-2,-1)(1.5,0){4}{\color{black!70}{\circle*{0.1}}}
\multiput(-2.75,-2)(1.5,0){4}{\color{black!70}{\circle*{0.1}}}
\multiput(-2,-3)(1.5,0){4}{\color{black!70}{\circle*{0.1}}}
\multiput(-2.75,-4)(1.5,0){4}{\color{black!70}{\circle*{0.1}}}
\multiput(-2,-5)(1.5,0){4}{\color{black!70}{\circle*{0.1}}}
\multiput(-1.25,-6)(1.5,0){3}{\color{black!70}{\circle*{0.1}}}
\multiput(-0.5,-7)(1.5,0){2}{\color{black!70}{\circle*{0.1}}}
\rput(0,-42){
\multiput(-1.25,-2)(1.5,0){2}{\color{black!70}{\circle*{0.1}}}
\multiput(-2,-1)(1.5,0){3}{\color{black!70}{\circle*{0.1}}}
\multiput(-1.25,0)(1.5,0){3}{\color{black!70}{\circle*{0.1}}}
\multiput(-2,1)(1.5,0){4}{\color{black!70}{\circle*{0.1}}}
\multiput(-2.75,2)(1.5,0){4}{\color{black!70}{\circle*{0.1}}}
\multiput(-2,3)(1.5,0){4}{\color{black!70}{\circle*{0.1}}}
\multiput(-2.75,4)(1.5,0){4}{\color{black!70}{\circle*{0.1}}}
\multiput(-2,5)(1.5,0){4}{\color{black!70}{\circle*{0.1}}}
\multiput(-1.25,6)(1.5,0){3}{\color{black!70}{\circle*{0.1}}}
\multiput(-0.5,7)(1.5,0){2}{\color{black!70}{\circle*{0.1}}}}
}}
\end{pspicture}  
\\
 $i\triangleleft j$ (resp  $i\triangleright j$). 
Obviously, $Q_0=\Lambda^{(1)}=\Lambda_{(n)}$  and $i\in \Lambda^{(j)}$ if and only if $j\in \Lambda_{(i)}$.
\\

 According to the Brauer-Humphreys reciprocity formulas $(P(j):\De(i))=[\nabla(i):S(j)]$ and $(I(j):\nabla(i))=[\De(i):S(j)]$ (see \cite{CPS})    the Axiom (2)  in the Definition~\ref{def1qh}  is equivalent to  the   analog   multiplicities axiom  for injective indecomposable  and costandard modules.
 For   any  1-quasi-hereditary algebra  $(A,\ma)$  and all $i,j\in Q_0(A)$  we thus  have 
\begin{center}
\hspace*{8mm}$(P(j):\De(i))=(I(j):\nabla(i))=[\De(i):S(j)]=[\nabla(i):S(j)]= \left\{
\begin{array}{cl}
	1 & \text{if } i\in \Lambda^{(j)},\\
	0 & \text{else.}
\end{array}
\right.$ \hspace{5mm} ($\ast$)
\end{center}

An  algebra   $\A$  is quasi-hereditary if and only if  the  opposite algebra $\A^{op}$ of $\A$   	related   to  the same partial order $\ma$ on  $Q_0(\A^{op})=Q_0(\A)$ 
 is  quasi-hereditary. 
 There    are  the   following   relationships between  the standard and costandard 	 as well as between   the $\De$-good and $\nabla$-good modules  of $\A$ and $\A^{op}$  
    (we denote  by $\DD $ the standard  $K$-duality): For all $i,j \in Q_0$,     we have   $\De_{\A}(i)\cong \DD(\nabla_{\A^{op}}(i))$  and  $[\De_{\A}(i):S(j)]=[\nabla_{\A^{op}}(i):S(j)]$.  
 For  $M\in \ger{F}(\De_{\A})$,  it is   $\DD(M)\in \ger{F}(\nabla_{\A^{op}})$ 
 and    $\left(M:\De(i)\right)=\left(\DD(M):\nabla_{\A^{op}}(i)\right)$.   
The corresponding dual properties  hold    for $\nabla_{\A}(i)$ and $M\in \ger{F}(\nabla_{\A})$. The  general properties of the standard duality imply  that  Axiom (3) and (4) in the Definition~\ref{def1qh} are self-dual   (see \cite[ Theorem 5.13]{ASS}). This   yields
  the following  lemma.

\begin{num}\begin{normalfont}\textbf{Lemma.}\end{normalfont} An algebra $A$ is 1-quasi-hereditary if  and only if   $A^{op}$  is   1-quasi-hereditary. 
\label{A0p}
\end{num}

\section{Projective indecomposables    and the $\Ext$-quiver}
\begin{small}
The 	structure of a 1-quasi-hereditary algebra $A$ is related to properties of the projective indecomposable modules, which will be  exhibited  in this section.  This implies  that the structure of the  standard $A$-modules, the quiver etc. is directly connected with the given  partial order. 
\end{small}
\\
  
The relationship between the  dimension vectors of an $A$-module  $M$  and of  the  subquotients  of  $M$ as well as  the equation $(\ast)$ shows  that  dimension  vectors of modules  $\De(j)$, $\nabla(j)$, $P(j)$, $I(j)$ and $A$  only  depend  on the structure of the  poset  $(Q_0,\ma)$.

 \begin{num}\begin{normalfont}\textbf{Lemma.}\end{normalfont} Let $A=(KQ/\II,\ma)$ be a 1-quasi-hereditary algebra and  $j,k\in Q_0$. Then 
\begin{itemize}
	\item[(1)] $\dimm_K\De(k)=\dimm_K\nabla(k)=\left|\Lambda_{(k)}\right|$,    \      \    $\displaystyle \dimm_KP(j)=\dimm_KI(j)=\sum_{k\in \Lambda^{(j)}}\left|\Lambda_{(k)}\right|$  \  and 
	\\
$\displaystyle \dimm_K A=\sum_{j\in Q_0}\left|\Lambda_{(j)}\right|^{2}$. 
	\item[(2)] $\displaystyle [P(j):S(k)]=[I(j):S(k)]= \left|\Lambda^{(j)}\cap \Lambda^{(k)}\right|$.
	\item[(3)] $P(1)\cong I(1)$, where  $1=\min \{Q_0,\ma\}$.
\end{itemize}
 \label{dimension} \end{num}
 
 \textit{Proof.} \textit{(1)} The dimensions   of  $\De(i)$, $\nabla(i)$, $P(i)$, $I(i)$ and  $A$ we obtain  directly from ($\ast$).

 \textit{(2)} The equation  ($\ast$) implies    $[P(j):S(k)]= \sum_{i\in \Lambda^{(j)}}[\De(i):S(k)]= \sum_{i\in \Lambda^{(j)}\cap \Lambda^{(k)}}[\De(i):S(k)] +  \sum_{i\in \Lambda^{(j)}\backslash \Lambda^{(k)}}[\De(i):S(k)]= \left|\Lambda^{(j)}\cap \Lambda^{(k)}\right|$. Similarly, we have  $[I(j):S(k)]=\left|\Lambda^{(j)}\cap \Lambda^{(k)}\right|$.
 
 \textit{(3)} The  Definition~\ref{def1qh} (3) implies   $P(1) \hookrightarrow I(1)$.  Since $\dimm_KP(1)\stackrel{\textit{(1)}}{=}\dimm_KI(1)$, we obtain $P(1)\cong I(1)$.
   \hfill $\Box$ 
\\

Any projective indecomposable  module over a 1-quasi-hereditary algebra  may be  considered as a submodule of $P(1)$ because of Definition~\ref{def1qh} (3) and Lemma ~\ref{dimension} \textit{(3)}.

\begin{num}\begin{normalfont}\textbf{Lemma.}\end{normalfont} Let $A=(KQ/\II,\ma)$ be a 1-quasi-hereditary algebra,   $i,j\in Q_0$  and $M(i)$  be a submodule of $P(1)$ isomorphic to  $P(i)$.   Then  
\begin{itemize}
	\item[(1)] $\Loc_{i} (M(j))\subseteq \Loc_{i}(M(i))$ 
	\item[(2)]  $\Loc_{i} (M(j))= \Loc_{i}(M(i))$  if and only if $i\in \Lambda^{(j)}$. 
\end{itemize}
In particular,  $P(i)\hookrightarrow P(j)$   if and only if   $i\in \Lambda^{(j)}$,  and there exists a  unique    submodule of $P(j)$ which  is    isomorphic to $P(i)$. 
  \label{lokale} \end{num}

\textit{Proof.} \textit{(1)}   Since    $\Loc_i (M(j)) = \left\{\left\langle m\right) \mid m \in M(j)_i\backslash\left\{0\right\}\right\}$  for all $i,j\in Q_0$,   it is
enough to show   $M(j)_i \subseteq M(i)_i$.
  Lemma~\ref{dimension} (2)  implies   $\dimm_KP(1)_i = \dimm_KM(i)_i = \left| \Lambda^{(i)}\right|$, thus  $M(i)\subseteq P(1)$ yields   $P(1)_i=M(i)_i$. Consequently,    $M(j)_i\subseteq P(1)_i=M(i)_i$ for all  $i,j\in Q_0$. 
  
 \textit{(2)} Obviously,  $i\in \Lambda^{(j)}$  if and only if $\left|\Lambda^{(i)}\cap \Lambda^{(j)}\right|=\left|\Lambda^{(i)}\right|$. In this case we have  $\dimm_KM(j)_i=\dimm_KM(i)_i$, thus   $M(j)_i \subseteq M(i)_i$  implies   $M(j)_i = M(i)_i$. 
 
Since  $\Loc_{i} (P(1))= \Loc_{i}(M(i))$,  we obtain that  for  any   submodule $N$ of $P(1)$ with $ N\cong P(i)$ it holds $N\subseteq M(i)$, thus  $\dimm_KN=\dimm_KM(i)$ implies   $N=M(i)$.  
Consequently,  $M(i)$ is the unique submodule  of $M(j)$ isomorphic to $P(i)$  if $i\in \Lambda^{(j)}$ because of    \textit{(2)}.
\hfill $\Box$

 \begin{num}\begin{normalfont}\textbf{Remark.} From now on,   for  $i,j\in Q_0$ with $i\in \Lambda^{(j)}$  we consider  $P(i)$ as a  submodule of $P(j)$. 
 Since  for every  $F\in \End_A(P(j))$ with  $F(P(i))\neq 0$  the submodule  $F(P(i))$ of $P(j)$ is  local with $\topp F(P(i))\cong S(i) $,   Lemma~\ref{lokale}  implies $F(P(i))\subseteq P(i)$. The submodule $P(i)$ of $P(j)$ is an $\End_A(P(j))^{op}$-module for all $i\in \Lambda^{(j)}$.
\end{normalfont}  
\label{p(i)}
    \end{num}

\begin{num}\begin{normalfont}\textbf{Lemma.}\end{normalfont}  Let $A=(KQ/\II,\ma)$ be a 1-quasi-hereditary algebra and  $j\in Q_0$. Then   
\begin{center}
$\displaystyle \De (j)= P(j)/\left(\sum_{j\triangleleft 
i} P(i)\right)$  \   \    and  \ \    $\displaystyle \nabla(j) = \bigcap_{j\triangleleft i} \Kern (I(j)\twoheadrightarrow I(i)).$
\end{center}
\label{stan}
    \end{num}
    \text{  }\\[-10mm]
    
    \textit{Proof.} Since   $\Loc_i(P(j))=\left\{\im (f) \mid f\in \Hom_A(P(i),P(j)), \ f\neq 0\right\}\stackrel{~\ref{lokale} }{=} \Loc_i(P(i))$,    we obtain   $\sum_{f\in \Hom_A(P(i),P(j))} \im (f) = P(i)$ for every $i\in \Lambda^{(j)}$. Moreover,  $\sum_{j<i}P(i)=\sum_{j\triangleleft i} P(i)$,  since  for every $k\in \Lambda^{(j)}\backslash \left\{j\right\}$  there exists $i\in Q_0$ with   $j\triangleleft i \ma k$,  thus   $P(k)\subseteq P(i)$. 
We obtain  $\De (j)\stackrel{~\ref{qhstandq}}{=} P(j)/\left(\sum_{j\triangleleft 
i} P(i)\right)$.  

Using the standard duality we have   $ \nabla(j) = \bigcap_{j\triangleleft i} \Kern (I(j)\twoheadrightarrow I(i)).$ \hfill $\Box$
\\

 Definition~\ref{def1qh} (4)  shows  that any standard module can be considered as a submodule of $\De(n)$.  Thus we consider  any submodule of $\De(j)$  as a submodule of $\De(n)$.

 \begin{num}\begin{normalfont}\textbf{Lemma.}\end{normalfont} Let    $A=(KQ/\II,\ma)$  be  a    1-quasi-hereditary, $j\in Q_0$. Then   $M$ is a submodule of $\De(j)$ if and only if  $\displaystyle M=\sum_{i\in \Lambda}\De(i)$ for some $\Lambda \subseteq \Lambda_{(j)}$. In particular,  $\Loc_i( \De(j)) =\left\{\De(i)\right\}$  if $i\in \Lambda_{(j)}$ and  $\Loc_i( \De(j)) =\emptyset$  if $i\not\in \Lambda_{(j)}$. Moreover, $\displaystyle \rad \De(j)=\sum_{j\triangleright i}\De(i)$.
\label{stan-inkl}
\end{num}
\text{  }\\[-10mm]

\textit{Proof.}    For every $i\in Q_0$ we have    $\Loc_i(\De(n))=\left\{\De(i)\right\}$, since $[\De(n):S(i)]=1$ (see Definition~\ref{def1qh} (2)).  If $i\in \Lambda_{(j)}$, then $[\De(j):S(i)]=1$, thus  $\Loc_i(\De(j))\neq \emptyset$. Since  $\Loc_i(\De(j)) \subseteq \Loc_i(\De(n))$,  we obtain   $\Loc_i(\De(j))=\left\{\De(i)\right\}$.  If $i\not\in \Lambda_{(j)}$, then $[\De(j):S(i)]=0$, thus  $\Loc_i(\De(j))= \emptyset$.
Any submodule $M$ of $\De(j)$ is a sum of some local submodules of $\De(j)$, thus $M=\sum_{i\in \Lambda}\De(i)$ for some $\Lambda \subseteq \Lambda_{(j)}$. In particularly,  $\rad \De(j) = \sum_{i\in \Lambda_{(j)}\backslash\left\{j\right\}} = \sum_{i\triangleleft j}\De(i)$, since for any $k\in \Lambda_{(j)}\backslash\left\{j\right\}$ there exists $i\in Q_0$ with $k\ma i\triangleleft j$, thus  $\De(k)\subseteq \De(i)$.
\hfill $\Box$ 
\\

\begin{num}\normalfont{\textbf{Remark.}}  Since  for a 1-quasi-hereditary algebra $A$ the algebra  $A^{op}$ is also 1-quasi-hereditary  (see ~\ref{A0p}),    every   statement  yields  a  corresponding  dual statement.     Lemma~\ref{stan-inkl} and   Lemma~\ref{lokale} implies that  for all  $j,l\in Q_0$  and all  $i\in \Lambda_{(j)}$ and $k\in \Lambda^{(l)}$ we obtain 
\\[5pt]
\hspace*{3mm}\begin{small}$S(1) \hookrightarrow \De(i)\hookrightarrow \De(j)\hookrightarrow P(k)\hookrightarrow P(l)\hookrightarrow P(1)\cong I(1)\twoheadrightarrow I(l) \twoheadrightarrow I(k)\twoheadrightarrow \nabla(j) \twoheadrightarrow \nabla(i) \twoheadrightarrow S(1). 
$ \end{small}
 \label{socel}
\end{num}

We are now going to determine the shape of the $\Ext$-quiver of a 1-quasi-hereditary algebra (cf. Theorem A (1)).

\begin{num}\begin{normalfont}\textbf{Theorem.} \end{normalfont}Let $A=\left(KQ/\II,\ma\right)$ be a 1-quasi-hereditary algebra. In the $\Ext$-quiver  of $A$ the vertices $i$ and $j$ are connected  by an arrow if and only if they are  neighbours with respect to $\ma$. Moreover,  if  $i\triangleleft j$ (or $i\triangleright j$)  then  $\left.\left|\left\{\alpha\in Q_1 \mid i \stackrel{\alpha}{\rightarrow}j \right\}\right| = \left|\left\{\alpha\in Q_1 \mid j \stackrel{\alpha}{\rightarrow}i \right\}\right|=1\right.$.
\label{topas} \end{num}

\text{  }\\[-15mm]

\textit{Proof.} Let $j,k\in Q_0$. The number of arrows from $k$ to  $j$   is the number of $S(k)$ in the decomposition of $\topp \left(\rad P(j)\right) $ (see  \cite[Lemma 2.12]{ASS}). We denote by   $\texttt{N}(j)$  the set $\left\{k\in Q_0\mid k\triangleleft j\right\} \cup \left\{k\in Q_0\mid  k\triangleright j \right\}$. 
We  have to show  $\topp \left(\rad P(j)\right)\cong 
\bigoplus_{k\in \texttt{N}(j)}S(k)$. In other
words, for every $k\in \texttt{N}(j)$  there exists $L(k)\in \Loc_k(P(j))$ with 
\begin{center}
$\rad P(j)=\sum_{k\in \texttt{N}(j)}L(k)$  \  \   and \  \   $L(t)\not\subseteq \sum_{k\in \texttt{N}(j)\atop t\neq k}L(k)$ \  \   for every $t\in \texttt{N}(j)$.
\end{center}
We denote   by   $\SM(\De(j))$ the set of  submodules of $\De(j)$    and by    $\SM(P(j)\mid \sum_{j\triangleleft i}P(i))$  the set of  submodules $M$ of $P(j)$ with  $\sum_{j\triangleleft i}P(i) \subseteq M$.
 The function  
$F: \SM\left(P(j)\mid \sum_{j\triangleleft 
i}P(i)\right) \rightarrow \SM(\De(j))$ with $F(M) = 
M/\left(\sum_{j\triangleleft i}P(i)\right)$
 is bijective (see ~\ref{stan}). By 
Lemma~\ref{stan-inkl}   for any  $k\in 
\Lambda_{(j)}$  there  exists $L(k)\in 
\Loc_k(P(j))$  such 
that $F\left(L(k)+\sum_{j\triangleleft i}P(i)\right) = \De(k)$       
  and  $F\left(\sum_{k\in \Lambda}L(k)+\sum_{j\triangleleft 
i}P(i)\right)=\sum_{k\in \Lambda}\De(k)$ for any  subset  $\Lambda \subseteq
\Lambda_{(j)}$, since  $F$ 
preserves and reflects inclusions.  
In particular,  $F\left(\rad P(j)\right)= \rad 
\De(j)=\sum_{j>k}\De(k)\stackrel{~\ref{stan-inkl}}{=}\sum_{j\triangleright k}\De(k)$, thus 
\begin{center}
$\rad 
P(j)=\sum_{j\triangleright k}L(k)+\sum_{j\triangleleft i}P(i)$. 
\end{center}
Since  $\De(t)\not\subseteq \sum_{j\triangleright k\atop t\neq k}\De(k)$,   we obtain    $L(t)\not\subseteq 
\sum_{j\triangleright k \atop t\neq k}L(k)+\sum_{j\triangleleft i}P(i)$ for every $t$ with  $j\triangleright t$. 	In order to prove   $P(t)\not\subseteq \sum_{j\triangleright k }L(k)+\sum_{j\triangleleft i \atop t\neq i}P(i)$ for $t$ with  $j\triangleleft t$,  it is enough to show the following two statements: Let $M,M'$ be some submodules of $P(j)$, then
\\[-3mm]

 \ding{192}  $P(t)\not\subseteq M$ and  $P(t)\not\subseteq M'$ implies  $P(t)\not\subseteq M+M'$,
 
 \ding{193}   $P(t)\not\subseteq L(k)$ for every  $k$ with  $j\triangleright k$ and $P(t)\not\subseteq P(i)$ for every  $i$ with  $j\triangleleft  i\neq t$.
 \\[5pt] 
\ding{192}: For all $m\in M_t$ and   $m'\in M'_t$ we have   $\left\langle m\right) \neq P(t)$ and   $\left\langle m'\right)\neq P(t)$. Since   $\left\langle m\right), \left\langle m'\right) \in  \Loc_{t}(P(j)) \stackrel{~\ref{lokale}}{\subseteq}\Loc_{t}(P(t))$  for $m,m'\neq 0$,   we obtain     $m, m'\in \rad P(t)$, thus   $m+m'\in \rad P(t)$. 
Consequently,  $\Loc_t(M+M')=\left\{\left\langle m+m'\right)\mid m\in M_t\backslash \left\{0\right\}, m'\in M_t'\backslash \left\{0\right\}\right\} \subseteq \Loc_t(\rad P(t))$ and   hence  $P(t)\not\subseteq M+M'$.
\\
\ding{193}: Assume $P(t)\subseteq L(k)$ for some $j\triangleright k$. Let  $G:P(k)\sur L(k)$, then   $L(k)\cong P(k)/\ker (G)$  implies the 	existence of $N\in \Loc_t(P(k))$ with $\ker (G)\subseteq N$ such that  $P(t)\cong N/\ker (G)$.  Since  $N\stackrel{ ~\ref{lokale}}{\subseteq} P(t)$, we have  $\ker G=0$ and $P(k)\cong L(k)$, a contradiction  because for   $j\triangleright k$, it holds  $P(k)\stackrel{~\ref{lokale}}{\not\subseteq} P(j)$. For all $i,t\in Q_0$ with   $j\triangleleft i,t$ and  $i\neq t$ we have $P(t)\not\subseteq P(i)$ by Lemma ~\ref{lokale}.
\hfill $\Box$

\section{A basis of a 1-quasi-hereditary algebra}
From now on  $A=(KQ/\II,\ma)$ is a 1-quasi-hereditary algebra  with $1\ma i\ma n$  for all 
$i\in Q_0$. We use the same notations as in the previous section.
\\[3pt]

The structure of the quiver of a 1-quasi-hereditary algebra shows that  for all $j,i,k\in Q_0$ with $i\in \Lambda^{(j)}\cap \Lambda^{(k)}$  there exists a path 
\begin{center}
 $
\begin{array}{l}
j\to \lambda_1 \to \cdots \to \lambda_m\to i$ \  with \  $j\ma \lambda_1\ma \cdots \ma \lambda_m\ma i \ \  \text{   resp.  }  \\
 i\to \mu_1 \to \cdots \to \mu_r\to k$ \  with \  $i\geqslant  \mu_1\geqslant \cdots \geqslant \mu_r\geqslant k
\end{array}$
\end{center}
 called \textit{increasing} path from $j$ to $i$,  resp. \textit{decreasing} path from $i$ to $j$. By concatenating these, we get a path from $j$ to $k$ passing through $i$, and we 
write $p(j,i,k)$  for the image in $A$ of such path. When 
$i=j=k$, the path $p(j,j,j)$ is the  trivial path $e_j$. 
All increasing resp.   decreasing 
paths  (as well all  arrows)  of the quiver occur in this way:  A  path of the form   $p(j,i,i)$   is  increasing resp.  $p(i,i,k)$ is a decreasing path.

 \begin{num} \normalfont{\textbf{Remark.}} Recall that the  module  generated by  $ p(j,i,k)$   is the image of the $A$-map  $f_{(j,i,k)}:P(k)\rightarrow P(j)$ via $f_{(j,i,k)}(e_k)=p(j,i,k)$, thus  a  submodule of $P(j)$  from $\Loc_k(P(j))$.

(a) Theorem~\ref{topas} implies   $ \rad P(j)= \sum_{j\triangleleft i}\left\langle j\to i\right) + \sum_{j\triangleright i}\left\langle j\to i\right)$ for any $j\in Q_0$.     Since $\left\langle j\to i\right)\in \Loc_i(P(j))$, we obtain that   $\left\langle j\to i\right)$ belongs to the submodule $P(i)$ of $P(j)$ for all $i$ with $j\triangleleft i$ (see ~\ref{lokale}).   It is easy to see that  $\left\langle j\to i\right)=P(i)$: 
 Assume $\left\langle j\to i\right) \subset P(i)$, then   $\left\langle j\to i\right)\not\subseteq \sum_{j\triangleleft i'\atop i\neq i'}\left\langle j\to i'\right) + \sum_{j\triangleright i}\left\langle j\to i\right)$  implies   $P(i)\not\subseteq \rad P(j)$ (see \ding{192} in the proof of ~\ref{topas}), a contradiction. 
 
  The $A$-map corresponding to  $(j\to i)$ with $j\triangleleft i$ is therefore  an inclusion. Consequently the $A$-map corresponding  to an increasing path $p(j,i,i)$  provides a composition of the  inclusions, thus  $f_{(j,i,i)}:P(i)\hookrightarrow P(j)$.  In particularly,  for any two increasing paths $p$ and $q$ from $j$ to $i$ we have $\left\langle p\right)=\left\langle q\right)$, since $\im (f_p)=\im (f_q) =P(i)$ (see ~\ref{lokale}).  Thus 
 $\left\langle p(j,i,i)\right)\cong \left\langle p(j',i,i)\right)\cong P(i)$ for all $j,j'\in \Lambda_{(i)}$.   Using  our notations,   we have  $\rad P(j)= \sum_{j\triangleleft i}P(i)+ \sum_{j\triangleright k} \left\langle p(j,j,k)\right) $.
  
(b) A path  $p(j,i,k) $ is the  product of   $p(i,i,k) $  and  $ p(j,i,i)$, therefore using (a) we have  
 $f_{(j,i,k)}:P(k) \stackrel{f_{(i,i,k)}}{\longrightarrow} P(i) \stackrel{f_{(j,i,i)}}{\hookrightarrow} P(j)$. Hence the    module $\left\langle p(j,i,k)\right)$  
 	 may  be considered as a  submodule of $P(i)(\subseteq P(j))$  from $\Loc_k(P(i))$.
 In particular,  it is easy to see that for all $j,k\in Q_0$ we have  $\left\langle p(j,n,k)\right)\cong \De(k)$ because $\De(k)$ is  the  uniquely    submodule of $P(n) = \De(n)$  from $\Loc_k(P(n))$ (see ~\ref{stan-inkl}).
   \label{projektiveiso} \end{num}

\begin{num} \begin{normalfont} \textbf{Theorem.}\end{normalfont} Let $A=(KQ/\II,\leqslant)$  be a 1-quasi-hereditary algebra and $j,k\in Q_0$. For any $i\in \Lambda^{(j)}\cap \Lambda^{(k)}$  we fix a path of the form $p(j,i,k)$.
   The set 
\begin{center}
$\left\{p(j,i,k) \mid i\in\Lambda^{(j)}\cap \Lambda^{(k)}\right\} \  $ is a $K$-basis of $P(j)_k$.
\end{center}
In particular,   $\ger{B}_j:= \left\{p(j,i,k) \mid  i\in\Lambda^{(j)}, \  k \in \Lambda_{(i)}\right\}   $ is a $K$-basis of $P(j)$ for any $j\in Q_0$ and $\ger{B}:=\left\{p(j,i,k) \mid j,k\in Q_0,  \ i\in\Lambda^{(j)}\cap \Lambda^{(k)}\right\}   $ is a $K$-basis of $A$. \label{theor}
\end{num}

\text{  }
\\
\parbox{13.2cm}{ The chosen  paths $p(j,i,k)$ for all $i\in \Lambda^{(j)}\cap \Lambda^{(k)}$ are symbolically represented in the picture to the  right (a path $p(j,i,k)$ is not uniquely  determined).  
 Theorem~\ref{theor}  shows that for any path 
 $p$ in $A$, which starts in $j$ and ends  in $k$  there exist $c_{i}\in K$ such that 
 $\displaystyle p=\sum_{i\in \Lambda^{(j)}\cap \Lambda^{(k)}}c_i\cdot p(j,i,k)$.  In other words  $\displaystyle p-\sum_{i\in \Lambda^{(j)}\cap \Lambda^{(k)}}c_i \cdot p(j,i,k) \in  \II $ (here $p$ and $p(j,i,k)$ are paths in  $Q$). 
 Using general methods   any   relation in $\II$   can by  transform  in this  form  (see for example \cite[Section II.2, II.3]{ASS}).  
}
\psset{xunit=0.67mm,yunit=0.72mm,runit=1mm}
\begin{pspicture}(0,0)(0,0)
\rput(28,6){
\begin{tiny}
\multiput(-1.25,5)(1.5,0){2}{\color{black!70}{\circle*{0.1}}}
\multiput(-2,4)(1.5,0){4}{\color{black!70}{\circle*{0.1}}}
\multiput(-4.25,3)(1.5,0){5}{\color{black!70}{\circle*{0.1}}}
\multiput(-5,2)(1.5,0){6}{\color{black!70}{\circle*{0.1}}}
\multiput(-5.75,1)(1.5,0){7}{\color{black!70}{\circle*{0.1}}}
\multiput(-5,0)(1.5,0){7}{\color{black!70}{\circle*{0.1}}}
\multiput(-5.75,-1)(1.5,0){7}{\color{black!70}{\circle*{0.1}}}
\multiput(-6.5,-2)(1.5,0){8}{\color{black!70}{\circle*{0.1}}}
\multiput(-5.75,-3)(1.5,0){7}{\color{black!70}{\circle*{0.1}}}
\multiput(-6.6,-4)(1.5,0){8}{\color{black!70}{\circle*{0.1}}}
\multiput(-5.75,-5)(1.5,0){7}{\color{black!70}{\circle*{0.1}}}
\multiput(-6.5,-6)(1.5,0){8}{\color{black!70}{\circle*{0.1}}}
\multiput(-5.75,-7)(1.5,0){7}{\color{black!70}{\circle*{0.1}}}
\multiput(-5,-8)(1.5,0){7}{\color{black!70}{\circle*{0.1}}}
\multiput(-4.25,-9)(1.5,0){6}{\color{black!70}{\circle*{0.1}}}
\multiput(-5,-10)(1.5,0){6}{\color{black!70}{\circle*{0.1}}}
\multiput(-4.25,-11)(1.5,0){6}{\color{black!70}{\circle*{0.1}}}
\multiput(-3.5,-12)(1.5,0){5}{\color{black!70}{\circle*{0.1}}}
\multiput(-2.75,-13)(1.5,0){4}{\color{black!70}{\circle*{0.1}}}
\multiput(-2,-14)(1.5,0){3}{\color{black!70}{\circle*{0.1}}}
\multiput(-1.25,-15)(1.5,0){3}{\color{black!70}{\circle*{0.1}}}
\begin{footnotesize}
\rput(0,10){\rnode{A}{$n$}}
\rput(0,-44){\rnode{B}{$1$}}
\rput(20,-17){\rnode{aa}{}}
\rput(-20,-17){\rnode{bb}{}}
\rput(-11,3){\rnode{A1}{}}
\rput(-11,-37){\rnode{B1}{}}
\rput(11,3){\rnode{A2}{}}
\rput(11,-37){\rnode{B2}{}}
\rput(13,16){\rnode{p0}{$p(j,i,k)$}}
\rput(-3,-9){\rnode{p1}{}}

              \rput(-7,-3){\rnode{0}{$i$}}
              \rput(3,-10){\rnode{01}{}}
              \rput(-10,-30){\rnode{1}{$j$}}
                        
              \rput(-9.4,-31.5){\rnode{j1}{}}
              \rput(-5,-35){\rnode{j2}{}}
              \rput(2,-32){\rnode{j23}{$p$}}
              \rput(-2,-32){\rnode{j3}{}}
              \rput(5,-37){\rnode{j4}{}}
              \rput(10,-25){\rnode{2}{$k$}}
              \rput(7,-10){\rnode{1r}{}}
              \rput(-11,-10){\rnode{2l}{}} \end{footnotesize}
              \end{tiny}
                      }
                      \ncline[linestyle=dashed]{-}{1}{j1}
                      \ncline[linestyle=dashed]{-}{j1}{j2}
                      \ncline[linestyle=dashed]{-}{j2}{j3}
                      \ncline[linestyle=dashed]{-}{j3}{j4}
                      \ncline[linestyle=dashed]{->}{j4}{2}
             \psset{nodesep=1pt,offset=0pt}
              \ncline{1}{0}
              \ncline{->}{0}{2}
              \psset{nodesep=1pt,offset=0pt}
              \ncline{1}{01}
              \ncline{->}{01}{2}
              \ncline{1}{A}
              \ncline{->}{A}{2}
              \psset{nodesep=1pt,offset=1pt}
              \ncline[linestyle=dotted]{21}{22}
              \ncline[linecolor=black!70]{->}{A}{A1}
              \ncline[linecolor=black!70]{->}{A1}{A}
              \ncline[linecolor=black!70]{->}{B1}{B}
              \ncline[linecolor=black!70]{->}{B}{B1}
              \ncline[linecolor=black!70]{->}{A}{A2}
              \ncline[linecolor=black!70]{->}{A2}{A}
              \ncline[linecolor=black!70]{->}{B2}{B}
              \ncline[linecolor=black!70]{->}{B}{B2}
              \ncarc[arcangle=-50,linestyle=dashed, linecolor=black!70]{A}{1}
              \ncarc[arcangle=25,linestyle=dashed, linecolor=black!70]{A}{1r}
              \ncarc[arcangle=25,linestyle=dashed, linecolor=black!70]{1r}{1}
              \ncarc[arcangle=35,linestyle=dashed, linecolor=black!70]{A}{2}
              \ncarc[arcangle=-35,linestyle=dashed, linecolor=black!70]{A}{2l}
              \ncarc[arcangle=-35,linestyle=dashed, linecolor=black!70]{2l}{2}
               \ncarc[arcangle=35,linecolor=black!70]{->}{p0}{p1}
                            
              \ncarc[arcangle=25,linestyle=dotted,linecolor=black!70]{A2}{aa}
              \ncarc[arcangle=-25,linestyle=dotted,linecolor=black!70]{->}{aa}{A2}
              \ncarc[arcangle=25,linestyle=dotted,linecolor=black!70]{->}{aa}{B2}
              \ncarc[arcangle=-25,linestyle=dotted,linecolor=black!70]{B2}{aa}
              \ncarc[arcangle=-25,linestyle=dotted,linecolor=black!70]{A1}{bb}
              \ncarc[arcangle=-25,linestyle=dotted,linecolor=black!70]{->}{bb}{B1}
              \ncarc[arcangle=25,linestyle=dotted,linecolor=black!70]{B1}{bb}
              \ncarc[arcangle=25,linestyle=dotted,linecolor=black!70]{->}{bb}{A1}
\end{pspicture} 
 \\[3pt]
   Part (2) of Theorem A follows directly from   ~\ref{theor}.
 \\

The proof of the Theorem ~\ref{theor} is based on the statements of the following  lemma. 
Recall that for  $i\in \Lambda^{(j)}$  we consider   $P(i)$ as a submodule of $P(j)$ (see ~\ref{p(i)}).

 \begin{num}\begin{normalfont}\textbf{Lemma.}\end{normalfont} Let $A=(KQ/\II,\leqslant)$  be a 1-quasi-hereditary algebra and $j,k\in Q_0$.  Let   $0\subset \cdots \subset D' \subset D \subset \cdots \subset P(j) $ be  a $\De$-good filtration of $P(j)$,   where   $D/D'\cong \De(i)$ for some $i\in \Lambda^{(j)} \cap \Lambda^{(k)}$. Then we have the following:  
\begin{itemize}
	\item[(1)] $D=P(i)+D'$.
		\item[(2)] $D'\subset \left\langle p(j,i,k)\right)+D' \subseteq D$  for any path of the form $p(j,i,k)$.
\end{itemize} 
In particular,  there exists a subset $\Lambda$ of $\Lambda^{(j)}$ with  $\displaystyle D=\sum_{i\in \Lambda}P(i)$. 
\label{filtrierung} \end{num}

\textit{Proof.} \textit{(1)} Let  $ 0=D(r+1) \subset D(r) \subset \cdots \subset D(1)=P(j) $ be a $\De$-good filtration  with   $D(l)/D(l+1)\cong \De(i_l)$ for all $r\leq l\leq 1$.  There is some local submodule $L(l)$ of $P(j)$ with top isomorphic to  $S(i_l)$ such that $D(l)= L(l)+D(l+1)$. Definition~\ref{def1qh} yields $ i_l\in \Lambda^{(j)}$  and therefore   $L(l)\subseteq P(i_l) \subseteq P(j)$   (see  ~\ref{lokale}).  We obtain   $D(l)=L(l)+D(l+1)\subseteq P(i_l)+D(l+1)$ for all $1\leq l\leq r$. 
In order to show  $D(l)= P(i_l)+D(l+1)$,   we have to show   $P(i_l)\subseteq D(l)$.

Assume $P(i_l)\not\subseteq D(l)$.  There exists $t\in \left\{1,\ldots, l-1\right\}$ with $P(i_l)\subseteq D(t)$  and $P(i_l)\not\subseteq D(t+1)$ and hence $D(t+1)\subset P(i_l)+ D(t+1) \subseteq D(t)$.  We show now   $P(i_l)+ D(t+1)=D(t)$, this  then implies   $D(t)/D(t+1)\cong \De(i_t)\cong \De(i_l)$ and hence  $\left(P(j):\De(i_l)\right)\geq 2$, a contradiction (see Definition~\ref{def1qh}).

  Since $0\neq  P(i_l)/\left(P(i_l)\cap  D(t+1)\right) \hookrightarrow D(t)/D(t+1)\cong \De(i_t)$,   the standard module $\De(i_t)$ has a local submodule with top isomorphic to  $S(i_l)$.  Thus $[\De(i_t):S(i_l)]\neq 0$ and hence $i_t\stackrel{~\ref{def1qh}}{\in} \Lambda^{(i_l)}$ and   therefore  $L(t)\subseteq P(i_t)\subseteq P(i_l)$ (see ~\ref{lokale}). Consequently,  $D(t)=L(t)+D(t+1)\subseteq P(i_l)+D(t+1) \subseteq D(t)$. We have $P(i_l)+D(t+1) = D(t)$.

 Via   induction on $r-k$ we obtain  $D(k)=\sum_{m=k}^{r}P(i_m)$ for any $0\leq k\leq r$. 
	
\textit{(2)} By  Lemma~\ref{stan}    and    \textit{(1)}, since $D/D'\cong P(i)/(P(i)\cap D')\cong \De(i)$,  we obtain  $P(i)\cap D'= \sum_{i \triangleleft l}P(l)$.   Because   $\left\langle p(j,i,k)\right)$  is a submodule of $P(i)\subseteq P(j)$ (see ~\ref{projektiveiso}(b)),    it  is    enough to  show $\left\langle p(j,i,k)\right) \not\subseteq    \sum_{i\triangleleft l}P(l)$. This implies     $\left\langle p(j,i,k)\right)\not\subseteq D'$   and  consequently     $D'\subset \left\langle p(j,i,k)\right)+ D' \subseteq P(i)+ D'=D$.

Let  $i\triangleright k$, then $p(i,i,k)=(i\to k)$.  
We have   $\left\langle p(i,i,k)\right) \not\subseteq \sum_{i\triangleleft  l}P(l)$, since    $\rad P(i)= \sum_{i\triangleright k}\left\langle p(i,i,k)\right)+\sum_{i\triangleleft  l}P(l)$ (see  ~\ref{projektiveiso}(a)). To deal with the general paths we consider   maps. 
Because  $\left\langle p(i,i,k)\right)\stackrel{~\ref{projektiveiso}}{=}\im (f_{(i,i,k)})$, we have
   $\im \left( P(k)\stackrel{f_{(i,i,k)}}{\rightarrow}P(i) \stackrel{\pi}{\sur} P(i)/\left(\sum_{i\triangleleft l}P(l)\right)\right) \neq 0$.   
 Since   $\im \left(\pi\circ f_{(i,i,k)}\right)\subseteq  P(i)/\left(\sum_{i\triangleleft l}P(l)\right)\stackrel{~\ref{stan}}{=}\De(i)$ and $\Loc_k(\De(i))\stackrel{~\ref{stan-inkl}}{=}\left\{\De(k)\right\}$,  we obtain   $\im \left(\pi \circ f_{(i,i,k)}\right) = \De(k)$. Lemma   ~\ref{stan}   implies   $\Kern \left(\pi \circ f_{(i,i,k)}\right) = \sum_{k\triangleleft j}P(j)$. This 	implies   a commutative diagram 
\begin{footnotesize}\begin{center}
$
\begin{array}{ccc}
	P(k) & \xrightarrow{f_{(i,i,k)}} & P(i) \\
	\downarrow & & \downarrow \pi\\
	\underbrace{P(k)/\left(\sum_{k\triangleleft j}P(j)\right)}_{\De(k)} & \xrightarrow{ \ \  \overline{f_{(i,i,k)}} \ \  } & \underbrace{P(i)/\left(\sum_{i\triangleleft l}P(l)\right)}_{\De(i)}
\end{array}
$
\end{center}
 \end{footnotesize}
The map   $\overline{f_{(i,i,k)}} $ is an inclusion,  since $\overline{f_{(i,i,k)}} \neq 0$.

Now  let  $i>k$ with $i\triangleright l_1\triangleright \cdots \triangleright l_m\triangleright k$. Inductively      we obtain  the  commutative diagrams  for the  path   $p(i,i,k)=\left(i\rightarrow l_1\rightarrow \cdots \rightarrow l_m\rightarrow k\right) = p(l_m,l_m,k)\cdot p(l_{m-1},l_{m-1},l_m)  \cdots  p(i,i,l_1)$ 
\begin{center}
\begin{footnotesize}
$
\begin{array}{cccccccccccc}
P(k)& \xrightarrow{f_{(l_m,l_m,k)}} & P(l_m) & \xrightarrow{f_{(l_{m-1},l_{m-1},l_m)}} &\cdots &        & \xrightarrow{f_{(l_1,l_1,l_2)}} & P(l_1) & \xrightarrow{f_{(i,i,l_1)}} & P(i) \\
\downarrow &      & \downarrow &         &   &    & &\downarrow & &\downarrow \pi \\
\Delta(k) & \hookrightarrow & \Delta(l_m) & \hookrightarrow & \cdots &  &      \hookrightarrow & \De(l_1)  & \hookrightarrow &  \Delta(i)
\end{array}
 $
\end{footnotesize}
\end{center}
For  the  maps  $f_{(i,i,k)}= f_{(i,i,l_1)} \circ f_{(l_1,l_1,l_2)} \circ \cdots  \circ  f_{(l_m,l_m,k)} $  and $\pi: P(i) \twoheadrightarrow P(i)/\left(\sum_{i\triangleleft i'}P(i')\right) \cong \De(i) $ we have $\im \left(\pi \circ f_{(i,i,k)}\right) \neq 0$, thus $\im(f_{(i,i,k)})=\left\langle p(i,i,k)\right)\not\subseteq  \sum_{i \triangleleft l}P(l)$. 
 Therefore   $f_{(j,i,k)}:P(k)\stackrel{f_{(i,i,k)}}{\longrightarrow} P(i)\stackrel{f_{(j,i,i)}}{\hookrightarrow}P(j)$ shows  that the submodule  $\im (f_{(j,i,k)}) =\left\langle p(j,i,k)\right)$  of $P(i)\subseteq P(j)$  is not the submodule of  $ \sum_{i \triangleleft l}P(l)$. \hfill $\Box$ 
\\

\textit{Proof of the theorem.} Let  $\ger{F}: 0=D(r+1) \subset D(r)\subset \cdots   \subset D(1)= P(j)$ be  $\De$-good, then 
$\left\{D(l)/D(l+1)\mid 1\leq l\leq r\right\}\stackrel{~\ref{def1qh}}{\longleftrightarrow}\left\{\De(i)\mid i\in \Lambda^{(j)}\right\}$.
Let   $\left\{i_1,\ldots , i_m\right\}=\Lambda^{(j)}\cap \Lambda^{(k)}$ such that   $\widetilde{\ger{F}} : 0\subseteq D(i_m+1)\subset D(i_m) \subseteq \cdots  \subseteq D(i_2+1)\subset D(i_2)\subseteq D(i_1+1)\subset D(i_1) \subseteq P(j)$  is a subfiltration of $\ger{F}$  with  
$D(i_t)/D(i_t+1)\cong \De(i_t)$ for   $1\leq t\leq m$.
By  Lemma~\ref{filtrierung} (2) the filtration $\widetilde{\ger{F}}$  can be refined to
\begin{center}
$
\begin{array}{cl}
	 0  & \subseteq D(i_m+1)\subset \left\langle p(j,i_m,k)\right)+  D(i_m+1) \subseteq D(i_m) \subseteq   \cdots\\
	 	&  \hspace{4cm}\vdots  \\
	  & \subseteq D(i_2+1)\subset \left\langle p(j,i_2,k)\right)+  D(i_2+1) \subseteq D(i_2)   \\
	 &  \subseteq D(i_1+1)\subset \left\langle p(j,i_1,k)\right)+  D(i_1+1) \subseteq D(i_1) \subseteq P(j)\\
\end{array}
$
\end{center}
Therefore    $p(j,i_1,k),\ldots , p(j,i_m,k)$ are linear independent in $P(j)_k$. Since  $m=$ $\left|\Lambda^{(j)}\cap \Lambda^{(k)}\right| $ $\stackrel{~\ref{dimension} (2)}{=} $ $\dimm_KP(j)_k $,       the  set $\left\{p(j,i,k) \mid i\in \Lambda^{(j)}\cap \Lambda^{(k)}\right\}  $ is a $K$-basis of $P(j)_k$.

Because  $\bigcup_{k\in Q_0}\left\{p(j,i,k) \mid i\in \Lambda^{(j)}\cap \Lambda^{(k)}\right\} = \left\{p(j,i,k)\mid i\in \Lambda^{(j)}, \ k\in \Lambda_{(i)}\right\}$, the set  $\ger{B}_j$ is a $K$-basis of $P(j)$. 
 \hfill $\Box$

\begin{num}\normalfont{\textbf{Remark.}} Let $j\in Q_0$ and $i,l\in \Lambda^{(j)}$ with $l\in \Lambda^{(i)}$, then  $p(j,l,k)\stackrel{~\ref{projektiveiso}(b)}{\in} P(l) \stackrel{~\ref{lokale}}{\subseteq}P(i)\stackrel{~\ref{lokale}}{\subseteq} P(j)$  for all $k\in \Lambda_{(l)}$.       We obtain that the set 
\begin{center}
$\ger{B}_j(i):= \left\{p(j,l,k)\mid l\in \Lambda^{(i)}, \ k\in \Lambda_{(l)}\right\}$  \ is a $K$-basis of the submodule $P(i)$ of $P(j)$,
\end{center}
since      $ \dimm_KP(i)\stackrel{~\ref{dimension}}{=} \sum_{l\in \Lambda^{(i)}} \left|\Lambda_{(l)}\right|=\left|\ger{B}_j(i)\right|$ and $\ger{B}_j(i)$ is a subset of $\ger{B}_j$ defined in  ~\ref{theor}. It is easy to check that 
 for all  subsets  $\Gamma_1, \Gamma_2$  of  $\Lambda^{(j)}$   and  $ \Gamma_{1,2}:=
  \left(\bigcup_{i\in \Gamma_1} \Lambda^{(i)}\right)\cap \left(\bigcup_{i\in \Gamma_2} \Lambda^{(i)}\right)$    the set 
 \begin{small}$ \left(\bigcup_{i\in \Gamma_1}\ger{B}_j(i)\right)\cap  \left(\bigcup_{i\in \Gamma_2}\ger{B}_j(i)\right) = \bigcup_{i\in\Gamma_{1,2}}\ger{B}_j(i) $\end{small}  is a $K$-basis of the submodule
\begin{center}
 $\left(\sum_{i\in \Gamma_1}P(i)\right)\cap \left(\sum_{i\in \Gamma_2}P(i)\right)= \sum_{i\in \Gamma_{1,2}}P(i)$ \ \ of  \ \ $P(j)$.
\end{center}
\label{basisproj} \end{num}

\section{Good filtrations}
\begin{small}In this section,  we  show   the relationship between    the   Jordan-Hölder filtrations of $\nabla(j)$ and    $\De$-good filtrations of $P(j)$ resp. the  Jordan-Hölder filtrations of $\De(j)$ and   $\nabla$-good filtrations of $I(j)$ over a 1-quasi-hereditary algebra $(A,\ma)$. The sets of these  Jordan-Hölder filtrations resp.    good filtrations are finite and related to  certain  sequences of elements from $\Lambda_{(j)}$  resp. $\Lambda^{(j)}$ which   depend  on $\ma$. 
\end{small}\\

For any $i\in \Lambda_{(j)}$ we can consider   the standard module $\De(i)$  as a submodule of $\De(j)$  and $\nabla(i)$ as a factor module of $\nabla(n)$ (see ~\ref{def1qh}(4)). 
 We denote by  $\ger{K}(j)$  the kernel of the   map $\nabla(n)\twoheadrightarrow \nabla(j)$. We have $\ger{K}(j)\subseteq \ger{K}(i)$ if and only if $i\in \Lambda_{(j)}$ (see ~\ref{socel}).  
 
 \begin{num}\begin{normalfont}\textbf{Proposition.}\end{normalfont} Let $A=\left(KQ/\II, \ma\right)$ be  1-quasi-hereditary, $j\in Q_0$,  $r=\left|\Lambda_{(j)}\right|$  and   $\mathcal{T}(j):=\left\{\textbf{i}=(i_1, i_2, \ldots , i_r) \mid i_m\in \Lambda_{(j)} , \   i_{k}\not\geqslant i_{t}, \ 1\leq  k < t \leq r \right\}$. Then   the following functions are bijective:
\begin{itemize}
\item[(1)] $ \mathcal{S}:\mathcal{T}(j)  \longrightarrow  \left\{\text{Jordan-Hölder-filtrations of } \De(j)\right\}$  with  
\\
\begin{small}$ \mathcal{S}(\textbf{i}): \ 0=J(0)\subset J(1)\subset  \cdots \subset J(t) \subset\cdots \subset J(r) \  \text{      such that  \   }   \displaystyle J(t):=\sum_{m=1}^{t}\De(i_m).
$ \end{small}
\\ Moreover,    $J(t)/J(t-1)\cong S(i_t)$  \  for  \    $1\leq t\leq r.$

\item[(2)] $\widetilde{\mathcal{S}}:\mathcal{T}(j)  \longrightarrow  \left\{\text{Jordan-Hölder-filtrations of } \nabla(j)\right\}$   with  
\\
\begin{small}$ \widetilde{\mathcal{S}}(\textbf{i}): \ \ger{J}(r)\subset   \cdots \subset \ger{J}(t) \subset\cdots \subset \ger{J}(1)\subset \ger{J}(0)=\nabla(j)   \text{   \   such that \    }  \displaystyle \ger{J}(t):= \left(\bigcap_{m=1}^{t}\ger{K}(i_m)\right)/\ger{K}(j).
$  \end{small}
\\ Moreover,   \  $\ger{J}(t-1)/\ger{J}(t)\cong S(i_t)$  \  for  \    $1\leq t\leq r.$
\end{itemize} \label{JoHo} \end{num}

\textit{Proof.} \textit{(1) }
By  definition of   $\mathcal{T}(j)$ for   $\textit{\textbf{i}}=(i_1,\ldots,i_t,\ldots,  i_r)\in \mathcal{T}(j)$ we have        $\Lambda_{(i_t)}\backslash \left\{i_t\right\} \subseteq \left\{i_1,\ldots , i_{t-1}\right\}$, 
 thus   $\rad \De(i_t) \stackrel{~\ref{stan-inkl}}{=}\sum_{l<i_t} \De(l) \subseteq J(t-1)$   for all $1\leq t\leq r$. For  $l\in \left\{i_1,\ldots , i_{t-1}\right\}$    we have  $[\De(l):S(i_t)]\stackrel{~\ref{def1qh}}{=}0$,  since $l\not\geqslant i_t$.  Because  $[\De(i_t):S(i_t)]=1$, we have   $\De(i_t)\not\subseteq J(t-1)$. Hence   $J(t)/J(t-1)\cong \De(i_t)/\left(\De(i_t)\cap J(t-1)\right)=\De(i_t)/\rad \De(i_t)\cong S(i_t) $  for every $1\leq t\leq r$.
 The function $ \mathcal{S}$  is well defined and injective.

Let $\mathcal{F}: \ 0=M(0)\subset M(1)\subset \cdots \subset M(r')=\De(j)$ be a  Jordan-Hölder-filtration of $\De(j)$   with $M(t)/M(t-1)\cong S(i_t)$ for all $1\leq t\leq r'$. Then $i_t\in \Lambda_{(j)}$ and   $r'\stackrel{~\ref{def1qh}}{=}\left|\Lambda_{(j)}\right|=r$.   There exists  $\Lambda(t)\subseteq \Lambda_{(j)}$  with     $M(t) \stackrel{~\ref{stan-inkl}}{=}\sum_{i\in \Lambda(t)}\De(i)$    for any $1\leq t\leq r$. By  induction  on $t$ we can show  $\Lambda(t)=\left\{i_1,\ldots , i_t\right\}$  with $i_k\not\geqslant i_v$ for $1\leq k < v \leq t$:
  Let $t=1$, then    $\De(i_1)=\De(1)\stackrel{~\ref{socel}}{=}\soc \De(j)$. Since   $M(t+1)/M(t)\cong S(i_{t+1})$,  we obtain  $\De(i_{t+1})\subseteq M(t+1)$ and $\De(i_{t+1})\not\subseteq M(t)$, because  
   $\Loc_{i_{t+1}}(\De(j))\stackrel{~\ref{stan-inkl}}{=}\left\{\De(i_{t+1})\right\}$. Thus $M(t+1)=M(t)+\De(i_{t+1})$ and  $l\not\geqslant i_{t+1}$ for all $l\in \left\{i_1,\ldots , i_t\right\}$.
   This implies  $\mathcal{F}=\mathcal{S}(i_1,\ldots, i_r)$, i.e.   the function  $\mathcal{S}$ is 	surjective.

 \textit{(2)} Since $A^{op}$ is also 1-quasi-hereditary (see ~\ref{A0p}),  by    duality  the function  $\mathcal{\widetilde{S}}$ is bijective. \hfill $\Box$ 
\\

In a similar way, we can  determine all   $\De$-good filtrations of  $P(j)$,  resp. $\nabla$-good filtrations of $I(j)$,  for every $j\in Q_0$.  For any $i\in \Lambda^{(j)}$ we continue denotind  by $P(i)$ the projective submodule of $P(j)$ with top isomorphic to $S(i)$ and by  $\mathcal{K}(j)$ we denote  the kernel of the map $P(1)\twoheadrightarrow I(j)$ (see ~\ref{socel}). Obviously,   it is      $\mathcal{K}(j)\subseteq \mathcal{K}(i)$ if and only if $i\in \Lambda^{(j)}$.

\begin{num}\begin{normalfont}\textbf{Proposition.}\end{normalfont} Let $A=\left(KQ/\II, \ma\right)$ be a 1-quasi-hereditary algebra, $j\in Q_0$,   $r=\left|\Lambda^{(j)}\right|$  and  $\mathcal{L}(j):=\left\{\textbf{i}=(i_1, i_2, \ldots , i_r) \mid i_m\in \Lambda^{(j)} , \   i_{k}\not\geqslant i_{t}, \ 1\leq  k < t \leq r \right\}$. Then   the following functions are bijective:
\begin{itemize}
	\item[(1)]
$
\mathscr{D}:\mathcal{L}(j)  \longrightarrow  \left\{\De\text{-good filtrations of } P(j)\right\}$  \   with  

\begin{small}$ 
\mathscr{D}(\textbf{i}) : \  0=D(r+1)\subset D(r)\subset \cdots \subset D(t) \subset \cdots \subset D(1) \   \text{  such that  }  \  \displaystyle D(t):=\sum_{m=t}^{r}P(i_m).
$  \end{small} \\ Moreover,   \  $D(t)/D(t+1)\cong \De(i_t)$ for every $1\leq t\leq r.$

\item[(2)] $
\mathscr{N}:\mathcal{L}(j)  \longrightarrow  \left\{\nabla\text{-good filtrations of } I(j)\right\}$  with 

$\mathscr{N}(\textbf{i}): \   N(1) \subset  \cdots \subset N(t) \subset\cdots \subset N(r) \subset I(j)  \  \text{with  }   \displaystyle N(t):=\left(\bigcap_{m=t}^{r}\mathcal{K}(i_m)\right)/\mathcal{K}(j).
$  \\  Moreover,  $N(t+1)/N(t)\cong \nabla(i_t)$ for every $1\leq t\leq r.$
\end{itemize} \label{degood}
\end{num}

\textit{Proof.} \textit{(1)} By  definition of   $\mathcal{L}(j)$ for   $\textit{\textbf{i}}=(i_1,\ldots , i_t , \ldots , i_r)\in \mathcal{T}(j)$ we have        $\Lambda^{(i_t)}\backslash \left\{i_t\right\} \subseteq \left\{i_{t+1},\ldots , i_r\right\}$ for any $1\leq t\leq r$. We obtain  $\ger{B}_j(i_t)\cap \left(\bigcup_{m=t+1}^{r}\ger{B}_j(i_m)\right)=\bigcup_{i_t<i}\ger{B}_j(i)$,  using the notations from ~\ref{basisproj}. Therefore  $P(i_{t})\cap D(t+1)= \sum_{i_{t}<i}P(i) =\sum_{i_{t}\triangleleft i}P(i)$ and consequently   $D(t)/D(t+1)\cong P(i_t)/\left(\sum_{i_{t}\triangleleft i}P(i)\right) \stackrel{~\ref{stan}}{=} \De(i_t)$. 
   The filtration $\mathscr{D}(\textit{\textbf{i}})$  is $\De$-good and  $\mathscr{D}$ is injective.

Let $\mathcal{F}: 0\subset D(r')\subset \cdots \subset D(1)=P(j)$ be a $\De$-good filtration, with $D(t)/D(t+1)\cong \De(i_t)$, then    $r'=r=\left|\Lambda^{(j)}\right|$  and    $D(t)=\sum_{m=t}^{r}P(i_m)$     (see Lemma~\ref{filtrierung}). 
The inclusion   $D(t) \subset D(k)$ implies  $P(i_{k})\not\subseteq P(i_t)$ for  $ k < t$. Hence         $i_k\not\geqslant i_t$  for all  $1\leq k<t\leq r$   (see  ~\ref{lokale})

\textit{(2)} follows from  the properies of  the standard duality   $\DD$. \hfill $\Box$ 
\\

The definitions of $\mathcal{T}(j)$ and $\mathcal{L}(j)$  yields $\mathcal{T}(n)=\mathcal{L}(1)$. Comparing  the  compositions factors of  the filtrations corresponding to $\textit{\textbf{i}}\in \mathcal{T}(n)$,    we obtain that  the  Jordan-Hölder filtration   $\widetilde{\mathcal{S}}(\textit{\textbf{i}})$ of $\nabla(n)$ induces the  $\De$-good  filtration  $\mathscr{D}(\textit{\textbf{i}})$ of $P(1)$.  Thus all $\De$-good filtrations of $P(1)$ can be represented in a diagram  whose shape   coincides with the submodule diagram of $\nabla(n)$. Moreover, any sequence from  $\mathcal{L}(j)$ can by completed to a sequence of $\mathcal{L}(1)$, thus  all $\De$-good filtrations of $P(j)$ are part  of this diagram for every $i\in Q_0$.  		Analogously,     the submodule diagram of $\De(n)$  and the diagram of all  $\nabla$-good filtrations of $I(1)$ has the same form  (for the illustration of  this  see Example 4 of  \cite{P}):
For  $\textit{\textbf{i}}=(i_1,\ldots , i_m,\ldots , i_n)\in \mathcal{T}(n)$  we have the following  	relationship     between  the  Jordan-Hölder filtrations of $\De(j)$ and $\De(n)$ (resp. $\nabla(j)$ and $\nabla(n)$)  as well as  $\De$-good filtrations of $P(j)$ and $P(1)$ (resp. $\nabla$-good filtrations of $I(j)$ and $I(1)$). The factors of the filtrations $ \mathcal{S}(\textit{\textbf{i}})$   and $\mathscr{N}(\textit{\textbf{i}})$ (resp. $\widetilde{\mathcal{S}}(\textit{\textbf{i}})$ and $\mathscr{D}(\textit{\textbf{i}})$) are labeled by the same vertices (as indicated above the corresponding filtration): 
\\
\begin{footnotesize}  
$
\begin{array}{lll}
 \mathcal{S}(\textit{\textbf{i}}): &\underbrace{\overbrace{0\subset J(1)}^{S(1)}\subset \cdots \subset \overbrace{J(t-1) \subset J(t)}^{S(i_t)}}_{\text{filtration of }  \De(j)}\subset  \cdots  \subset \overbrace{J(n-1)\subset J(n)}^{S(n)}=\De(n) & \text{ for } \ (i_1,\ldots , i_t)\in \mathcal{T}(j),\\[20pt]
\widetilde{\mathcal{S}}(\textit{\textbf{i}}):	 & \overbrace{\ger{J}(n)\subset \ger{J}(n-1)}^{S(n)}\subset \cdots \subset \underbrace{\overbrace{\ger{J}(t)\subset \ger{J}(t-1)}^{S(i_t)}\subset \cdots \subset \overbrace{\ger{J}(1)\subset \ger{J}(0)}^{S(1)}}_{\text{filtration of } \nabla(j)}=\nabla(n) & \text{ for }  \  (i_1,\ldots , i_t)\in \mathcal{T}(j),\\[20pt]
\mathscr{D}(\textit{\textbf{i}}):	&  \underbrace{\overbrace{0\subset D(n)}^{\De(n)}\subset \cdots \subset \overbrace{D(t+1)\subset D(t)}^{\De(i_t)}}_{\De\text{-good  filtration of }  P(j)}\subset \cdots \subset \overbrace{D(2)\subset D(1)}^{\De(1)}=P(1)   & \text{ for  }  \  (i_t,\ldots , i_n)\in \mathcal{L}(j),\\[20pt]
\mathscr{N}(\textit{\textbf{i}})	& \overbrace{ N(1)\subset N(2)}^{\nabla(1)}\subset \cdots \subset \underbrace{\overbrace{N(t)\subset N(t+1)}^{\nabla(i_t)}\subset \cdots \subset \overbrace{N(n)\subset N(n+1)}^{\nabla(n)}}_{\nabla\text{-good  filtration of }  I(j)}=I(1)   & \text{ for  } \ (i_t,\ldots , i_n)\in \mathcal{L}(j),
\end{array}
$
\end{footnotesize}  
\\

Let  $\Lambda \subseteq Q_0$ and   $\check{\Lambda}:=\bigcup_{i\in \Lambda}\Lambda^{(i)}$. We can always  construct a sequence $(i_1,\ldots , i_t,\ldots , i_n)\in \mathcal{L}(1)$     with     $\left\{i_t,\ldots , i_n\right\}=\check{\Lambda}$. For any $k\in \check{\Lambda}$  there exists an $i\in \Lambda$ with $i\ma k$, thus  $P(k)\subseteq P(i)$ and consequently 
 $\sum_{l\in \Lambda}P(l)=\sum_{l\in \check{\Lambda}}P(j)$.

\begin{num}\begin{normalfont}\textbf{Corollary.}\end{normalfont} Let $\Lambda_1$  and $\Lambda_2$  be some subsets of $Q_0$  with $\check{\Lambda}_2\subset \check{\Lambda}_1$. Then for  the submodules $\displaystyle M_1:= \sum_{l\in \Lambda_1}P(l)$ and $\displaystyle M_2:=\sum_{l\in \Lambda_2}P(l)$     of $P(1)$, it is $M_2\subset M_1$  and $M_1/M_2\in \ger{F}(\De)$ (resp. for  $\displaystyle N_1:= \bigcap_{l\in \Lambda_1} \mathcal{K}(l) $  and $\displaystyle  N_2:= \bigcap_{l\in \Lambda_2}\mathcal{K}(l)$ we have  $N_1\subset N_2$   and $N_2/N_1\in \ger{F}(\nabla)$)  with 
\begin{center}
$\displaystyle \left(M_1/M_2:\De(k)\right)=\left\{
\begin{array}{cl}
	1 & \text{if  } k\in \check{\Lambda}_1\backslash \check{\Lambda}_2, \\
	0 & \text{else}
\end{array}
\right.$    \ \  and  \  \   $\displaystyle \left(N_2/N_1:\nabla(k)\right)=\left\{
\begin{array}{cl}
	1 & \text{if  } k\in \check{\Lambda}_1\backslash \check{\Lambda}_2, \\
	0 & \text{else.}
\end{array}
\right.$
\end{center}
\label{multfactor} \end{num}

\textit{Proof.} 
  We can construct  a  sequence $\textit{\textbf{i}}=(i_1, \ldots , i_{t_1},\ldots , i_{t_2}, \ldots , i_n)\in \mathcal{L}(1)$ such that $\left\{i_{t_v},\ldots , i_n\right\}=\check{\Lambda}_v$ for $v=1,2$.  In the $\De$-good filtration $\mathscr{D}(\textit{\textbf{i}})$  of $P(1)$   we have $D(t_v)=\sum_{m=t_v}^{n}P(i_l)=M_v$  and a $\De$-good filtrations of $M_v$ for $v=1,2$ 
  \begin{center}
$\mathscr{D}(\textit{\textbf{i}}): \underbrace{\underbrace{0\subset D(n) \subset \cdots \subset D(t_2)}_{\De\text{-good filtration of } M_2}\subset  \cdots \subset D(t_1)}_{\De\text{-good filtration of } M_1}\subset \cdots\subset D(1)=P(1).$
\end{center}
 Since  $\check{\Lambda}_2\subset \check{\Lambda}_1$,   we have  $M_2\subset M_1$ and  $M_1/M_2\in \ger{F}(\De)$ because  the induced filtration  $D(t_2)/M_2\subset D(t_2-1)/M_2\subset \cdots \subset D(t_1)/M_2$  is $\De$-good.   
 The properties of the filtration $\mathscr{D}(\textit{\textbf{i}})$ implies   $\left(M_v:\De(l)\right)=1$  for all $l\in \check{\Lambda}_v$ and $\left(M_v:\De(l)\right)=0$ for all $l\in Q_0\backslash\check{\Lambda}_v$, here $v=1,2$. Thus $\left(M_1/M_2:\De(k)\right)=\left(M_1:\De(k)\right)-\left(M_2:\De(k)\right)$  implies the statement. 
 
  The dual statement  follows by  dual argumentation.  \hfill $\Box$
\\

For every quasi-hereditary algebra $\A$,   the  category $\ger{F}(\De)$     is a \textit{resolving} subcategory of $\modd \A$ (resp. $\ger{F}(\nabla)$  is a  \textit{coresolving}  subcategory of $\modd \A$), i.e. the  category $\ger{F}(\De)$
 is closed under extensions,   kernels of surjective maps and it  contains all projective $\A$-modules (resp. $\ger{F}(\nabla)$ is closed under extensions,  cokernels  of injective  maps and contains all injective  $\A$-modules)   (see   \cite[Theorem 3 (resp. Theorem  3*)]{Rin1}). 
 
  Using  this fact,  when dealing with 1-quasi-hereditary algebras   we can determine  all local modules in $\ger{F}(\De)$ resp. colocal  modules in $\ger{F}(\nabla)$.

\begin{num}\begin{normalfont}\textbf{Corollary.}\end{normalfont} Let $A=(KQ/\II,\ma)$ be a 1-quasi-hereditary algebra,  $j\in Q_0$ and   $M,N$ be $A$-modules with $\topp M\cong S(j)$, $\soc N \cong S(j)$.  Then 
\begin{itemize}
	\item[(1)] 
$M\in \ger{F}(\De)$  if and only if $\displaystyle M\cong P(j)/\left(\sum_{i\in \Lambda}P(i)\right)$ for some $\Lambda\subseteq \Lambda^{(j)}\backslash \left\{j\right\}$.

	\item[(2)] 
$N\in \ger{F}(\nabla)$  if and only if $\displaystyle N\cong \left(\bigcap_{i\in \Lambda}\mathcal{K}(i)\right)/\mathcal{K}(j)$ for some $\Lambda\subseteq\Lambda^{(j)}\backslash \left\{j\right\}$.  
\end{itemize}
\label{gutfiltr}
\end{num}

\textit{Proof.} \textit{(1)} The filtration $0\subseteq  \ker\left(P(j)\twoheadrightarrow M\right) \subset P(j) $ can be refined to  a $\De$-good filtration  $\mathscr{D}(\textit{\textbf{i}})$  for some  $\textit{\textbf{i}}\in \mathcal{L}(j)$, thus $\ker\left(P(j)\twoheadrightarrow M\right) =\sum_{i\in \Lambda}P(i)$ for some $\Lambda\subseteq \Lambda^{(j)}$. 
Since  $M\neq 0$, we have $j\not\in \Lambda$. The other direction follows from  Corollary~\ref{multfactor}.

\textit{(2)}  is the dual statement of \textit{(1)}.  \hfill $\Box$

\begin{num}\normalfont{\textbf{Remark.}}  If  for all $j\in Q_0$ and all $i\in \Lambda^{(j)}$,   and  any two paths $p$, $q$   of the form  $p(j,i,i)$  it is  $p=q$ and any two paths $p'$, $q'$  of the form $p(i,i,j)$ it is  $p'=q'$,   
 then the   algebra    $B(A)=KQ_{B(A)}/I_{B(A)} \ $  given by the quiver  $Q_{B(A)}=(Q_0, \left\{(j\to i)\in Q_1 \mid j<i\right\})$   with all commutativity relations 
 and the partially ordered set $(Q_0,\leqslant)$  
 is a (strong) 
 exact Borel subalgebra   of $A$  and  $C(A):=B(A)^{op}$ is a  $\De$-subalgebra of $A$   in the sense of König (see [K]). 
The structure of  the  $A$-module   $\De(j)$  corresponds to the  structure  of   $P_{C(A)}(j)$  (this also holds for   $\nabla(j)$   and  $I_{B(A)}(j)$).
In this case ~\ref{degood} is  a  consequence  of \cite[Proposition 2.5]{K}  and ~\ref{JoHo}.  

All  known  1-quasi-hereditary  algebras   have   exact  Borel and $\De$-subalgebras. We  conjecture that this is in general the case.
\end{num}

\section{The characteristic tilting module }
\begin{small}
For  any quasi-hereditary basic algebra $(\A,\ma)$  the full subcategory  $\ger{F}(\De)\cap \ger{F}(\nabla)$ of $\modd \A$   consisting of all $\A$-modules which are $\De$-good and $\nabla$-good is  determined    by the  so called \textit{characteristic tilting module}  $T_{\A}$ of $A$ defined by Ringel in \cite{Rin1}:   For any $i\in Q_0$ there exists  an   (up to isomorphism) uniquely determined indecomposable $\A$-module $T_{\A}(i)$ in  $\ger{F}(\De)\cap \ger{F}(\nabla)$ with the following properties: For $j\not\ma i$ it is   $(T_{\A}(i):\De(j)) =(T_{\A}(i): \nabla (j))= [T_{\A}(i):S(j)]=0$ and  $(T_{\A}(i):\De(i))=(T_{\A}(i):\nabla(i)) =[T_{\A}(i):S(i)] =1$. Moreover,  there exists   a submodule $Y_{\A}(i)\in \ger{F}(\nabla)$  of $T_{\A}(i)$ with $T_{\A}(i)/Y_{\A}(i)\cong \nabla(i)$   
 (resp.  a factor module $X_{\A}(i)\in \ger{F}(\De)$  with  $\Kern (T_{\A}(i)\twoheadrightarrow X_{\A}(i)) \cong \De(i)$).   The $A$-module  $T_{\A} $  is isomorphic to  $ \bigoplus_{i\in Q_0}T_{\A}(i)$.   Moreover,  any module in $\ger{F}(\De)\cap \ger{F}(\nabla)$ is a  direct sum of some copies  of $T_{\A}(i)$.

   We recall the notations and properties  of  some  factor algebra  of a quasi-hereditary  algebra $\A=(KQ/\II,\ma)$, which will be used later: Let  $\Lambda$ be some saturated  subset of $Q_0$ (i.e. if $v\in \Lambda$ and $k\in Q_0$ then $k<v$ implies $k\in \Lambda$),  by $J(\Lambda)$ we denote   the ideal   $\A(\sum_{i\in Q_0\backslash \Lambda}e_i)\A$  of $\A$. For the  quiver $Q(\Lambda)$  of the  factor algebra $\A(\Lambda) := \A/J(\Lambda)$ we have    $Q_0(\Lambda)=\Lambda$  and 
$Q_1(\Lambda)=\left\{(i\rightarrow j)\in Q_1\mid i,j\in \Lambda\right\}$. 
 All paths $p=\left(k_1\rightarrow k_2\rightarrow \cdots \rightarrow k_m\right)$ in  $\A$  with  $k_t \not\in  \Lambda$ for some $1\leq t\leq m$   span   $J(\Lambda)$ as a $K$-space. 
  Moreover,  all    $\A(\Lambda)$-modules can be   considered as  the $\A$-modules  $M$ with $[M:S(i)]=0$ for all $i\in Q_0\backslash \Lambda$. The  projective   $\A(\Lambda)$-module $P_{\A(\Lambda)}(i)$ is isomorphic to the $\A$-module  $P(i)/J(\Lambda)P(i)$ for every $i\in Q_0(\Lambda)$. 
In particular,  the algebra $(\A(\Lambda),\ma)$  is quasi-hereditary 
 with $\De(i)\cong \De_{\A(\Lambda)}(i)$ and $\nabla(i)\cong \nabla_{\A(\Lambda)}(i)$ for all $i\in \Lambda$  (see   \cite{DD}). We have  $\ger{F}(\De_{\A(\Lambda)})\subseteq \ger{F}(\De)$ (resp. $\ger{F}(\nabla_{\A(\Lambda)}) \subseteq \ger{F}(\nabla)$)   and $T_{\A(\Lambda)}$ is a direct summand of $T_{\A}$  (more precisely $T_{\A(\Lambda)}(i)\cong T_{\A}(i)$).
 \end{small}
\\

Let $A=(KQ/\II,\ma)$ be a 1-quasi-hereditary algebra with $1\ma i\ma n$ for all $i\in Q_0$.  
Since    $P(1)\stackrel{~\ref{dimension}}{\cong} I(1)$   admits   $\De$-good and  $\nabla$-good filtrations with  $X(n)= P(1)/P(n) \in \ger{F}(\De)$ and   $Y(n)= \Kern \left(P(1)\twoheadrightarrow I(n)\right) \in \ger{F}(\nabla)$ (see ~\ref{gutfiltr}), we have $P(1)\cong T(n)$.  
\\[3pt]
\parbox{12cm}{We fix  $i\in Q_0$.  The factor algebra $A(i)$  of $A$ is defined  as follows:

\begin{center}
$\displaystyle A(i):= A/J(i)  \   $ where  $ \  \displaystyle J(i):= A\left(\sum_{j\in Q_0\backslash \Lambda_{(i)}}e_j\right)A$.
\end{center}
}
\psset{xunit=1.1mm,yunit=1.3mm,runit=3.5mm}
\begin{pspicture}(18,0)(0,0)
\rput(22,9){
\begin{tiny}
\rput(0,6){\rnode{nn}{$n$}}
\rput(-1,5.6){\rnode{0n}{}}
\rput(1,5.6){\rnode{0n0}{}}
\rput(-10.5,-6){\rnode{0nl}{}}
\rput(10.5,-6){\rnode{0nr}{}}
              \rput(0,0.5){\rnode{0i}{$i$}}
              \rput(0,0){\rnode{0}{\textcolor{white}{...\ .\ ,,,.}}}
              \rput(-5,-3){\rnode{1}{$i_1$}}
              \rput(5,-3){\rnode{2}{$i_{t'}$}}
              \rput(-3,-6){\rnode{3}{$i_1$}}
              \rput(0,-6){\rnode{33}{$\cdots $}}
              \rput(3,-6){\rnode{4}{$i_r$}}
              \rput(15,-5.8){\rnode{55}{$ Q(i)  $}}
              \rput(13,-6.2){\rnode{qui}{}}
              \rput(8.3,-7.8){\rnode{qui2}{}}
              \rput(0,-9){\rnode{5}{$j$}}
              \rput(-3,-12){\rnode{6}{$k_1$}}
              \rput(0,-12){\rnode{66}{$\cdots $}}
              \rput(3,-12){\rnode{7}{$k_{m}$}}
              \rput(0,-18){\rnode{00}{\textcolor{white}{..\ .\ .,,, .}}}
              \rput(0,-18.5){\rnode{00c}{$1$}}
              \rput(-5,-15){\rnode{01}{$j_{1}$}}
              \rput(5,-15){\rnode{02}{$j_{t}$}}      \end{tiny}     }
              \psset{nodesep=1.3pt,offset=1.3pt,linewidth=0.5pt,arrows=<-}
              \ncarc[arcangle=-35,linestyle=dashed]{1}{01}
              \ncarc[arcangle=35]{qui2}{qui}
              \ncarc[arcangle=35,linestyle=dashed]{01}{1}
              \ncarc[arcangle=-25,linestyle=dashed]{0}{3}
              \ncarc[arcangle=25,linestyle=dashed]{3}{0}
              \ncarc[arcangle=-25,linestyle=dashed]{6}{00}
              \ncarc[arcangle=25,linestyle=dashed]{00}{6}
              
              \psset{offset=1.3pt,nodesep=0pt,arrows=<-}
              \ncarc[arcangle=39,linestyle=dashed,linecolor=black!80]{00}{0nl}
              \ncarc[arcangle=39,linestyle=dashed,linecolor=black!80]{-}{0nl}{0n}
              \ncarc[arcangle=-39,linestyle=dashed,linecolor=black!80]{-}{0nl}{00}
              \ncarc[arcangle=-39,linestyle=dashed,linecolor=black!80]{0n}{0nl}
              
              \ncarc[arcangle=-39,linestyle=dashed,linecolor=black!80]{00}{0nr}
              \ncarc[arcangle=-39,linestyle=dashed,linecolor=black!80]{-}{0nr}{0n0}
              \ncarc[arcangle=39,linestyle=dashed,linecolor=black!80]{-}{0nr}{00}
              \ncarc[arcangle=39,linestyle=dashed,linecolor=black!80]{0n0}{0nr}
              \ncline{0}{1}
              \ncline{1}{0}
               \ncline{3}{5}
               \ncline{5}{3}
               \ncline{5}{6}
               \ncline{6}{5}
               \ncline{5}{7}
               \ncline{7}{5}
               \ncline{00}{01}
               \ncline{01}{00}
              \psset{nodesep=1.3pt,offset=1.3pt,linewidth=0.5pt,arrows=->}
              \ncarc[arcangle=35,linestyle=dashed]{2}{02}
              \ncarc[arcangle=-35,linestyle=dashed]{02}{2}
              \ncarc[arcangle=-25,linestyle=dashed]{4}{0}
              \ncarc[arcangle=25,linestyle=dashed]{0}{4}
              \ncarc[arcangle=25,linestyle=dashed]{7}{00}
              \ncarc[arcangle=-25,linestyle=dashed]{00}{7}
              \ncline{0}{2}
              \ncline{2}{0}
              \ncline{4}{5}
              \ncline{5}{4}
              \ncline{00}{02}
               \ncline{02}{00}
              
\end{pspicture}\\[3pt]
For the  quiver $Q(i)$ of $A(i)$ we have    $Q_0(i):=\Lambda_{(i)}$ and   $Q_1(i):=\left\{(j\to k)\in Q_1 \mid j,k\in \Lambda_{(i)}\right\}$.

\begin{num}\begin{normalfont}\textbf{Theorem}.\end{normalfont} \textit{Let $A=(KQ/\II, \ma)$ be a   1-quasi-hereditary algebra   and  $i \in Q_0$.   The following statements are equivalent:
\begin{itemize}
	\item[(i)] $A(i)$ is 1-quasi-hereditary,
	\item[(ii)] $\displaystyle T(i)\cong P(1)/\left(\sum_{l \in Q_0\backslash \Lambda_{(i)}}P(l)\right)$, \hspace{1.3cm}
	(ii') $\displaystyle T(i)\cong \bigcap_{{l \in Q_0\backslash \Lambda_{(i)}}}\Kern (I(1)\twoheadrightarrow
 I(l)),$
 \item[(iii)] $\displaystyle \soc T(i)$ is simple, \hspace{3.45cm}  
(iii') 	$\displaystyle\topp T(i)$ is simple.
\end{itemize}} \label{A(i)} \end{num}

The subset $\Lambda_{(i)}$ of $Q_0$ is saturated, thus $(A(i),\ma)$ is a quasi-hereditary algebra. 
The proof of this theorem is based on some properties of projective $A(i)$-modules, which  we consider  in the next lemma. 
For   $A(i)$-modules resp. paths  we   use   the index $(i)$. It should be noted that for any $l\in \Lambda_{(i)}$ a path  $p(j,l,k)$  runs through some vertices from $\Lambda_{(i)}$ (see Sec.3).

\begin{num}\begin{normalfont}\textbf{Lemma.}\end{normalfont}  Let    $i\in Q_0$ and    $(A(i),\ma)$ be 	 defined as above.  Then the following statements  hold for any $j\in  Q_0(i)$.
\begin{itemize}
\item[(a)] $\displaystyle P_{(i)}(j) \cong P(j)/\left(\sum_{l\in \Lambda^{(j)} \backslash  \Lambda_{(i)}} P(l)\right)$  and 
\\[5pt]
$\left\{p_{(i)}(j,l,k) \mid l\in \Lambda^{(j)} \cap \Lambda_{(i)} , \ k\in \Lambda_{(l)}\right\}$  is a $K$-basis of $P_{(i)}(j),$  \psset{xunit=0.38mm,yunit=0.38mm,runit=1mm}
\begin{pspicture}(0,0)(0,0)
\rput(58,58){
\rput(0,0){\begin{tiny}
\rput(0,10){\rnode{A}{$n$}}
\rput(0,-44){\rnode{B}{$1$}}
\rput(17,-17){\rnode{aa}{}}
\rput(-17,-17){\rnode{bb}{}}
              \rput(0,-7){\rnode{01}{$i$}}
             \begin{footnotesize} \rput(45,5){\rnode{lin}{$\Lambda^{(j)}\backslash \Lambda_{(i)}$}} \end{footnotesize}
              \rput(25,4){\rnode{ll1}{}}
              \rput(6,-1){\rnode{ll2}{}}
              \rput(0,-26){\rnode{2}{$j$}}
              \rput(8,-10){\rnode{1r}{}}
              \rput(-12,-10){\rnode{2l}{}}
              \end{tiny}
                             \ncarc[arcangle=-35]{->}{ll1}{ll2}
        \ncarc[arcangle=65,linestyle=dashed]{A}{2}
        \ncarc[arcangle=-65,linestyle=dashed]{A}{2}
          \ncarc[arcangle=65,linestyle=dashed]{01}{B}
        \ncarc[arcangle=-65,linestyle=dashed]{01}{B}
              \ncline[linestyle=dotted]{21}{22}
              \ncline[linecolor=black]{-}{A}{A1}
              \ncline[linecolor=black]{-}{A1}{A}
              \ncline[linecolor=black]{-}{B1}{B}
              \ncline[linecolor=black]{-}{B}{B1}
              \ncline[linecolor=black]{-}{A}{A2}
              \ncline[linecolor=black]{-}{A2}{A}
              \ncline[linecolor=black]{-}{B2}{B}
              \ncline[linecolor=black]{-}{B}{B2}
              \ncarc[arcangle=35,linestyle=dotted,linecolor=black]{A}{aa}
              \ncarc[arcangle=-35,linestyle=dotted,linecolor=black]{-}{aa}{A}
              \ncarc[arcangle=35,linestyle=dotted,linecolor=black]{-}{aa}{B}
              \ncarc[arcangle=-35,linestyle=dotted,linecolor=black]{B}{aa}
              \ncarc[arcangle=-35,linestyle=dotted,linecolor=black]{A}{bb}
              \ncarc[arcangle=35,linestyle=dotted,linecolor=black]{-}{bb}{A}
              \ncarc[arcangle=-35,linestyle=dotted,linecolor=black]{-}{bb}{B}
              \ncarc[arcangle=35,linestyle=dotted,linecolor=black]{B}{bb}
\rput(0,-3){
\multiput(-2,3)(1.5,0){3}{\color{black!80}{\circle*{0.1}}}
\multiput(-1.25,2)(1.5,0){3}{\color{black!80}{\circle*{0.1}}}
\multiput(-2,1)(1.5,0){4}{\color{black!80}{\circle*{0.1}}}
\multiput(-2.75,0)(1.5,0){5}{\color{black!80}{\circle*{0.1}}}
\multiput(-3.5,-1)(1.5,0){5}{\color{black!80}{\circle*{0.1}}}
\multiput(-2.75,-2)(1.5,0){5}{\color{black!80}{\circle*{0.1}}}
\multiput(-2,-3)(1.5,0){2}{\color{black!80}{\circle*{0.1}}}
\multiput(-2.75,-4)(1.5,0){1}{\color{black!80}{\circle*{0.1}}}
\multiput(-3.5,-5)(1.5,0){1}{\color{black!80}{\circle*{0.1}}} 
\multiput(2.5,-3)(-1.5,0){1}{\color{black!80}{\circle*{0.1}}} 
\multiput(3.25,-4)(-1.5,0){1}{\color{black!80}{\circle*{0.1}}}
 }
}}
\end{pspicture}
\psset{xunit=0.38mm,yunit=0.38mm,runit=1mm}
\begin{pspicture}(0,0)(0,0)
\rput(55,-5){
\rput(0,0){\begin{tiny}
\rput(0,10){\rnode{A}{$n$}}
\rput(0,-44){\rnode{B}{$1$}}
\rput(17,-17){\rnode{aa}{}}
\rput(-17,-17){\rnode{bb}{}}
              \rput(0,-7){\rnode{01}{$i$}}
             \begin{footnotesize} \rput(47,-10){\rnode{lin}{$\Lambda^{(j)}\cap \Lambda_{(i)}$}} \end{footnotesize}
              \rput(25,-11){\rnode{ll1}{}}
              \rput(6,-16){\rnode{ll2}{}}
              \rput(0,-26){\rnode{2}{$j$}}
              \rput(8,-10){\rnode{1r}{}}
              \rput(-12,-10){\rnode{2l}{}}
              \end{tiny}
               \ncarc[arcangle=-35]{->}{ll1}{ll2}
        \ncarc[arcangle=65,linestyle=dashed]{A}{2}
        \ncarc[arcangle=-65,linestyle=dashed]{A}{2}
        \ncarc[arcangle=65,linestyle=dashed]{01}{B}
        \ncarc[arcangle=-65,linestyle=dashed]{01}{B}
              \ncline[linestyle=dotted]{21}{22}
              \ncline[linecolor=black]{-}{A}{A1}
              \ncline[linecolor=black]{-}{A1}{A}
              \ncline[linecolor=black]{-}{B1}{B}
              \ncline[linecolor=black]{-}{B}{B1}
              \ncline[linecolor=black]{-}{A}{A2}
              \ncline[linecolor=black]{-}{A2}{A}
              \ncline[linecolor=black]{-}{B2}{B}
              \ncline[linecolor=black]{-}{B}{B2}
              \ncarc[arcangle=35,linestyle=dotted,linecolor=black]{A}{aa}
              \ncarc[arcangle=-35,linestyle=dotted,linecolor=black]{-}{aa}{A}
              \ncarc[arcangle=35,linestyle=dotted,linecolor=black]{-}{aa}{B}
              \ncarc[arcangle=-35,linestyle=dotted,linecolor=black]{B}{aa}
              \ncarc[arcangle=-35,linestyle=dotted,linecolor=black]{A}{bb}
              \ncarc[arcangle=35,linestyle=dotted,linecolor=black]{-}{bb}{A}
              \ncarc[arcangle=-35,linestyle=dotted,linecolor=black]{-}{bb}{B}
              \ncarc[arcangle=35,linestyle=dotted,linecolor=black]{B}{bb}
\rput(0,-16){
\multiput(-0.5,5)(1.5,0){2}{\color{black!80}{\circle*{0.1}}}
\multiput(-1.25,4)(1.5,0){3}{\color{black!80}{\circle*{0.1}}}
\multiput(-2,3)(1.5,0){4}{\color{black!80}{\circle*{0.1}}}
\multiput(-2.75,2)(1.5,0){5}{\color{black!80}{\circle*{0.1}}}
\multiput(-2,1)(1.5,0){4}{\color{black!80}{\circle*{0.1}}}
\multiput(-1.25,0)(1.5,0){3}{\color{black!80}{\circle*{0.1}}} 
\multiput(-0.5,-1)(1.5,0){2}{\color{black!80}{\circle*{0.1}}} }
}}
\end{pspicture}
\item[(b)] $P_{(i)}(j) \hookrightarrow P_{(i)}(1)$,   \   $I_{(i)}(1)\sur I_{(i)}(j)$,
\item[(c)] $\left(P_{(i)}(j):\De_{(i)}(k)\right)= [\De_{(i)}(k):S_{(i)}(j)]=1 $ for all $k\in \Lambda^{(j)}\cap\Lambda_{(i)}$.
\end{itemize} \label{factoralgebra} \end{num}

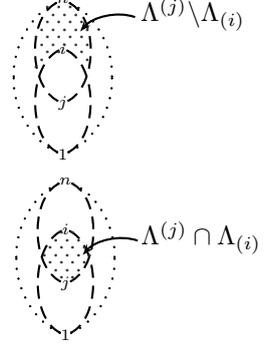
	\textit{Proof}.\textit{(a)}  Since  \begin{small}$P_{(i)}(j)\cong P(j)/(J(i)P(j))$,\end{small}  it is   enough to show \begin{small}$J(i)P(j)= \sum_{l\in \Lambda^{(j)}\backslash \Lambda_{(i)}} P(l)$.\end{small} 
 The set  $\underbrace{\left\{p(j,l,k) \mid l\in \Lambda^{(j)} \backslash \Lambda_{(i)} , \ k\in \Lambda_{(l)}\right\}}_{\textbf{B}_1:=}\cup \underbrace{\left\{p(j,l,k) \mid l\in \Lambda^{(j)} \cap \Lambda_{(i)} , \ k\in \Lambda_{(l)}\right\}}_{\textbf{B}_2:=}$ is a $K$-basis of $P(j)$ (see  ~\ref{theor}).		
Any  path   starting  in $j$ and  passing through some $l\in Q_0\backslash \Lambda_{(i)}$  belongs to $\spann_K\textbf{B}_1$. Thus $\textbf{B}_1=\bigcup_{l\in \Lambda^{(j)} \backslash \Lambda_{(i)} }\ger{B}_j(l)$ is a $K$-basis of $J(i)P(j)$ and of  the submodule $\sum_{l\in \Lambda^{(j)} \backslash \Lambda_{(i)}}P(l)$ of $P(j)$, in the notation  of ~\ref{basisproj}.
 We have $J(i)P(j)= \sum_{l\in \Lambda^{(j)}\backslash \Lambda_{(i)}} P(l)$  and 
	  $\left\{p_{(i)}(j,l,k) \mid p(j,l,k)\in \textbf{B}_2 \right\}= \left\{p_{(i)}(j,l,k) \mid l\in \Lambda^{(j)} \cap \Lambda_{(i)} , \ k\in \Lambda_{(l)}\right\}$ is a $K$-basis of $P_{(i)}(j)$.
	 
 \textit{(b)}  We have  $P(j)\cap \left(\sum_{l\in Q_0\backslash \Lambda_{(i)}}P(l)\right)= \sum_{l\in \Lambda^{(j)}\backslash \Lambda_{(i)}}P(l)$, according  to ~\ref{basisproj}  for the subsets  $\Gamma_1=\Lambda^{(j)}$ and $\Gamma_2=Q_0\backslash\Lambda_{(i)}$ of $Q_0=\Lambda^{(1)}$. Thus  
\\
\begin{small}$\underbrace{\displaystyle P(j)/\left(\sum_{l\in \Lambda^{(j)}\backslash \Lambda_{(i)}}P(l)\right)}_{P_{(i)}(j)}\cong \displaystyle \left(P(j)+\sum_{l\in Q_0\backslash \Lambda_{(i)}}P(l)\right)/\left(\sum_{l\in Q_0\backslash \Lambda_{(i)}}P(l)\right) \subseteq \underbrace{P(1)/\left(\sum_{l\in Q_0\backslash \Lambda_{(i)}}P(l)\right)}_{P_{(i)}(1)}$. \end{small}
\\
 Therefore  $P_{(i)}(j)$ can be  considered as a  submodule of $P_{(i)}(1)$ for any $j\in \Lambda_{(i)}$.

Any  projective indecomposable $A(i)^{op}$-module can by embedded in the projective indecomposable $A(i)^{op}$-module corresponding to the minimal vertex $1$ because $A^{op}$ is 1-quasi-hereditary and $A^{op}(i)\cong A(i)^{op}$. Using duality, we obtain $I_{(i)}(1)\sur I_{(i)}(j)$.

\textit{(c)} Since $\De(k)\cong \De_{(i)}(k)$,  we have $[\De_{(i)}(k):S_{(i)}(j)]\stackrel{~\ref{def1qh}}{=}1$ for all $k\in \Lambda^{(j)}\cap\Lambda_{(i)}$.
 For the sets $\Lambda_1=\Lambda^{(j)}=\check{\Lambda}_1$ and $\Lambda_2=\Lambda^{(j)}\backslash \Lambda_{(i)}=\check{\Lambda}_2$  (in the notation of ~\ref{multfactor}) 
we have   $M_1=P(j)$ and $M_2=\sum_{l\in \Lambda^{(j)}\backslash \Lambda_{(i)}}P(l)$. Thus   $\left(M_1/M_2:\De_{(i)}(k)\right)=1$ for all $k\in \check{\Lambda}_1\backslash \check{\Lambda}_2= \Lambda^{(j)}\cap\Lambda_{(i)}$.
\hfill $\Box$ 
\\

For all  $j\in Q_0(i)$ it is   $1\ma j\ma i$  and  $\De_{(i)}(j)\cong \De(j)$ as well as $\nabla_{(i)}(j)\cong \nabla(j) $, thus  $\De_{(i)}(j) \hookrightarrow \De_{(i)}(i)$ and  $\nabla_{(i)}(i)\sur \nabla_{(i)}(j) $  (see ~\ref{socel}).
The  foregoing   lemma shows that  the axioms of a 1-quasi-hereditary algebra are satisfied  for  $(A(i),\ma)$  if and only if  if and only if $P_{(i)}(1)\cong I_{(i)}(1)$.   
\\

\textit{Proof of the  theorem.}  Let $i\in Q_0$.  Since $\soc \De(j)\stackrel{~\ref{socel}}{\cong} S(1)$  for all $i\in Q_0$ and $T(i)\in \ger{F}(\De)$,  we obtain $\soc T(i)\cong S(1)^{m}$ for some $m\geq 1$.

$(i)\Rightarrow (ii)$ If $A(i)$ is 1-quasi-hereditary, then  $ P_{(i)}(1)\cong I_{(i)} (1)$   is isomorphic to  $T_{(i)}(i) $, since     $i$   is   maximal in $Q_0(i)$. 
The $A$-modules     $T_{(i)}(i)$ and $T(i)$ are isomorphic.  
  Lemma~\ref{factoralgebra} \textit{(a)}  implies  $T(i)\cong P_{(i)}(1)\cong P(1)/\left(\sum_{l\in Q_0 \backslash  \Lambda_{(i)}} P(l)\right)$.
 
 $(ii)\Rightarrow (iii)$     Since      $\soc T(i)\cong S(1)^{m}\cong \nabla(1)^m$  and  $ T(i)\in \ger{F}(\nabla)$, the filtration $0\subset \soc T(i) \subset T(i)$ can be refined  to a $\nabla$-good filtration of $T(i)$  since $\ger{F}(\nabla)$ is coresolving.
 We have  $\big(T(i):\nabla(1)\big) = \big(\soc T(i):\nabla(1)\big)+\big(T(i)/\soc T(i):\nabla(1)\big)$.    It is enough to show $\big(T(i):\nabla(1)\big)=1$, this implies   $[\soc T(i):S(1)]=\left(\soc T(i):\nabla(1)\right)=m=1$.
 
  Since     $T(i)\cong P(1)/\left(\sum_{l\in Q_0 \backslash  \Lambda_{(i)}} P(l)\right) $, \  $P(1)\in  \ger{F}(\De)$,       the filtration $ 0\subset \sum_{l\in Q_0 \backslash  \Lambda_{(i)}} P(l) \subset P(1)$ can be refined to a $\De$-good filtration  $\mathscr{D}(\textit{\textbf{i}})$    for some $\textit{\textbf{i}}=(i_1,\ldots ,i_t ,\ldots , i_n) \in  \mathcal{L}(1)$ (see ~\ref{degood}).  There  exists  $1\leq t< n$ with $D(t+1)=\sum_{l\in Q_0 \backslash  \Lambda_{(i)}} P(l)$. Thus  $\left(T(i):\De(j)\right)=1$ for $j\in \left\{i_1,\ldots , i_t\right\}$ and $\left(T(i):\De(j)\right)=0$ for $j\in \left\{i_{t-1},\ldots , i_n\right\}$. In the notation of ~\ref{multfactor}  for $\Lambda_1=Q_0$ and $\Lambda_2=Q_0\backslash \Lambda_{(i)}$  we obtain $T(i)\cong M_1/M_2$
and  
$(T(i):\De(j))=\begin{small}\left\{
\begin{array}{cl}
	1 & \text{if } j\in \Lambda_{(i)}, \\
	0 & \text{else}.
\end{array}
\right. \end{small}$ 
 Hence $\Lambda_{(i)}=\left\{i_1,\ldots, i_t\right\}$.  Since  $i_k\not\geqslant i_v$ for  $1\leq k < v \leq t$,  we have  $(i_1,\ldots , i_t)\in \mathcal{T}(i)$ (see ~\ref{JoHo}). Thus $i_1 =1$ and  $i_t=i$.
   
    Let now  $\mathscr{N}: 0=N(r+1)\subset N(r)\subset  \cdots \subset   N(1)=T(i)$ be a $\nabla$-good filtration  with   $N(v)/N(v+1)\cong \nabla(j_v)$ for every $1\leq v\leq r$. We have to  show $\left\{i_1,\ldots, i_t\right\} \subseteq \left\{j_1,\ldots ,j_r\right\}$. Then  the filtrations $\mathscr{D}(\textit{\textbf{i}})$ and  $\mathscr{N}$ as well as $\dimm_K\De(j)\stackrel{~\ref{dimension}}{=}\dimm_K\nabla(j)$   implies  
\begin{center}
$\displaystyle \dimm T(i)=\sum_{j\in \left\{i_1,\ldots, i_t\right\}}\dimm_K\De(j)= \sum_{j\in \left\{i_1,\ldots, i_t\right\}}\dimm_K\nabla(j) + \underbrace{\sum_{j\in \left\{j_1,\ldots, j_r\right\} \backslash \left\{i_1,\ldots, i_t\right\} }\dimm_K\nabla(j)}_{=0}. $ 
\end{center}
In other words, this implies  $\left\{i_1,\ldots, i_t\right\} = \left\{j_1,\ldots ,j_r\right\}$ and $t=r$. Consequently,   for all $j\in \left\{i_1,\ldots, i_t\right\} $ we obtain  $\left(T(i):\nabla(j)\right)=1$ and therefore  $\left(T(i):\nabla(i_1)\right)=\left(T(i):\nabla(1)\right)=1.$
 
 We show this  by induction on $t-w$: If $w=0$, then $i=i_t\in \left\{j_1,\ldots ,j_r\right\}$, since $\left(T(i):\nabla(i)\right)=1$  by the properties of $T(i)$. Assume $i_{t-w},i_{t-(w-1)},\ldots , i_t\in \left\{j_1,\ldots ,j_r\right\}$.
   For the $k$-th coordinate  of the dimension vector of  $T(i)$   we have 
\\[5pt]
$\begin{array}{ccll}
\displaystyle [T(i):S(k)] & = & \displaystyle  \sum_{l\in \{i_{1},\ldots , i_{t-(w+1)}\}}[\De(l):S(k)] + \sum_{j\in \{i_{t-w},\ldots , i_t\}}[\De(j):S(k)] &  \begin{scriptsize}(\De\text{-good filtration } \mathscr{D}(\textit{\textbf{i}}))\end{scriptsize}\\[10mm]
& = & \displaystyle  \sum_{j\in \{j_1,\ldots , j_r\} \atop j\not\in \{i_{t-w},\ldots , i_t\}}[\nabla(j):S(k)] +  \sum_{j\in \{i_{t-w},\ldots , i_t\}}[\nabla(j):S(k)]   &  \begin{scriptsize} (\nabla\text{-good filtration } \mathscr{N})\end{scriptsize}\\
\end{array}$ 
\\[5pt]
Let $X(k):=\sum_{l\in \{i_{1},\ldots , i_{t-(w+1)}\}}[\De(l):S(k)]$ and $Y(k):=\sum_{j\in \{j_1,\ldots , j_r\} \atop j\not\in \{i_{t-w},\ldots , i_t\}}[\nabla(j):S(k)]$ for $k\in Q_0$. 
Since    $[\De(j):S(k)]=[\nabla(j):S(k)]$ for all $j,k\in Q_0$ (see Sec.1 $(\ast)$), we  obtain $X(k)=Y(k)$ for all $k\in Q_0$. By  definition of $\mathcal{T}(i)$  for  $(i_1,\ldots , i_{t-(w+1)}, \ldots ,i_t )\in \mathcal{T}(i)$  we obtain   $ i_1,\ldots , i_{t-(w+1)} \not\in  \Lambda^{(i_{t-(w+1)})}\backslash \left\{i_{t-(w+1)}\right\}$. Thus $X(k)\stackrel{~\ref{def1qh}}{=}0=Y(k)$ for all $k\in \Lambda^{(i_{t-(w+1)})}\backslash \left\{i_{t-(w+1)}\right\}$. We obtain $\{j_1,\ldots , j_r\} \backslash \{i_{t-w},\ldots , i_t\}\not\subseteq \Lambda^{(i_{t-(w+1)})}\backslash \left\{i_{t-(w+1)}\right\}$. Moreover,   for  $k=i_{t-(w+1)}$ we have $X(k)\neq 0$    since $[\De(k):S(k)]=1$, therefore   $Y(k)\neq 0$. There exists $j\in \{j_1,\ldots , j_r\} \backslash \{i_{t-w},\ldots , i_t\}$ with   $[\nabla(j):S(i_{t-(w+1)})]=1$, hence   $j \in \Lambda^{(i_{t-(w+1)})}$.  Thus $j\not\in \Lambda^{(i_{t-(w+1)})}\backslash \left\{i_{t-(w+1)}\right\}$  and $j\in \Lambda^{(i_{t-(w+1)})} $ implies  $j=i_{t-(w+1)}\in \left\{j_1,\ldots , j_r\right\}$.

 $(iii)\Leftrightarrow (ii')$ The socle of $ T(i)$ is simple  if and only if $T(i)$ is a submodule of $I(1)$.  
 The filtration $0\subset T(i) \subset I(1)$ can be refined to a $\nabla$-good filtration $\mathscr{N}(\textit{\textbf{i}})$  for some $\textit{\textbf{i}}=(i_1,\ldots ,i_t,\ldots ,  i_n)\in \mathcal{L}(1)$ (see ~\ref{degood})  
 There exists     $1\leq t\leq n$  with  $T(i)=
 N(t) = \bigcap_{m=t}^{n}\Kern (I(1)\twoheadrightarrow I(i_m))$  and    $(T(i):\nabla(j))=\begin{small} \left\{
\begin{array}{cl}
1 & \text{if } j\in \left\{i_{1},\ldots , i_{t-1}\right\}, \\
0 & \text{if } j\in \left\{i_{t},\ldots , i_{n}\right\}.
\end{array}
\right.\end{small}$  
  We know that    $T(i)$    satisfies   $\left(T(i):\nabla(j)\right)=0$ for all $j\in Q_0\backslash \Lambda_{(i)}$ and   $(T(i):\nabla(j))\neq 0$ implies $j\in \Lambda_{(i)}$.  Since $i_1,\ldots , i_{t-1}\in \Lambda_{(i)}$ and $i\not\in \left\{i_t,\ldots , i_n\right\}$,  we obtain  $\Lambda_{(i)}\cap \left\{i_{t},\ldots , i_n\right\} =\emptyset $. Therefore  $\left\{i_{t},\ldots , i_n\right\}=Q_0\backslash \Lambda_{(i)}$.

 $(ii')\Rightarrow (iii')$ The dual argumentation of $(ii)\Rightarrow (iii)$. 
 
  $(iii')\Leftrightarrow (ii)$ The dual argumentation of $(iii)\Leftrightarrow (ii')$.
  
  $(iii)\Rightarrow (i)$ If $\soc T(i)\cong S(1)$, then  $(iii)\Rightarrow (ii')\Rightarrow (iii')\Rightarrow (ii)$  implies   $T(i)\cong P(1)/\left(\sum_{l\in Q_0 \backslash  \Lambda_{(i)}} P(l)\right) \stackrel{~\ref{factoralgebra} (a)}{\cong} P_{(i)}(1)$.  Since  $\soc P_{(i)}(1)\cong S(1)$ and $\dimm_KP_{(i)}(1)=\dimm_KI_{(i)}(1)$    (see Brauer-Humphreys reciprocity formulas and Lemma~\ref{factoralgebra} \textit{(c)}), we obtain  $P_{(i)}(1)\cong I_{(i)}(1)$.   Therefore   
   the algebra   $A(i)$  is 1-quasi-hereditary.\hfill $\Box$

\begin{num}\normalfont{\textbf{Remark.}} 
If $i\in Q_0$ is a neighbor of 1 (i.e. $1\triangleleft i$), then  for the    $A(i)$-module  $P_{(i)}(1)$  we have   $\rad P_{(i)}(1)=P_{(i)}(i)\cong \De_{(i)}(i) \cong \De(i)$ because $0\subset P_{(i)}(i)\subset P_{(i)}(1)$ is the uniquely  determined $\De$-good filtration of $P_{(i)}(1)$. Therefore $\soc P_{(i)}(1)\cong S(1) $ and  consequently $A(i)$ is 1-quasi-hereditary.   Theorem~\ref{A(i)} implies  that for any 1-quasi-hereditary algebra  $A=(KQ/\II, \ma)$ with $1\ma i\ma n$ it is:
\begin{itemize}
	\item $T(1)\cong \De(1)\cong \nabla(1)\cong S(1)$, 
	\item $T(n)\cong P(1)\cong I(1)$,
	\item $\displaystyle T(i) \cong P(1)/\left(\sum_{j\in Q_0\backslash \left\{1,i\right\}}P(j)\right)\cong \bigcap_{j\in Q_0\backslash \left\{1,i\right\}}\Kern(P(1)\twoheadrightarrow I(i))$ for any $i\in Q_0$ with $1\triangleleft i$.
\end{itemize}
  \label{tilting} \end{num}
An example of a 1-quasi-hereditary algebra $A$ such that  for some $i\in Q_0(A)$ the algebra $A(i)$ is not 1-quasi-hereditary can be found in \cite{P}.

\section{The Ringel dual of a  1-quasi-hereditary algebra} 
\begin{small}
The concept of Ringel duality 
    is  specific   to the theory of quasi-hereditary algebras (see \cite{Rin1}): For any quasi-hereditary (basic) algebra $\A$ the endomorphism algebra of the  characteristic tilting $\A$-module $T_{\A}$ is called the \textit{Ringel dual} of $\A$,  denoted by $R(\A)$  [i.e. $R(\A)=\End_{\A}\left(T_{\A}\right)^{op}$].  Since the direct summands of $T_{\A}$ are pairwise   non isomorphic, $R(\A)$ is a basic algebra. The vertices in  the quiver $Q(R(\A))$ may be identified with the vertices  of  $Q(\A)$ [$T(i) \leftrightsquigarrow i$]. The algebra $R(\A)$  is quasi-hereditary with respect to the opposite order on $Q_0(\A)$.  Furthermore,   $R(R(\A))$ and   $\A$ are  isomorphic as  quasi-hereditary algebras. The functor
 $\mathscr{R}_{(\A)}(-):=\Hom_{\A}(T_{\A}, -): \modd \A \rightarrow  \modd R(\A)$ induces an equivalence between $\ger{F}_{\A}(\nabla)$ and $\ger{F}_{R(\A)}(\De)$ 
and    for any $i\in Q_0(\A)$ hold
\begin{center}
$\mathscr{R}_{(\A)}(\nabla_{\A}(i))=\Delta_{R(\A)}(i), \  \   \    $  $   \mathscr{R}_{(\A)}(T_{\A}(i))=P_{R(\A)}(i), \   \  \    $ $ \mathscr{R}_{(\A)}(I_{\A}(i))=T_{R(\A)}(i).$
\end{center}
Applying $\mathscr{R}_{(\A)}(-)$ to an exact sequence  $0\to M'\to M\to M''\to 0$ in $\modd \A$ with $M',M,M''\in \ger{F}(\nabla)$ yields an exact sequence in $R(\A)$   
 and $\left(M:\nabla(i)\right)=\left(\mathscr{R}_{(\A)}(M):\De_{R(\A)}(i)\right)$ for all $i\in Q_0(\A)$.
 \end{small}
\\

The next theorem shows that the class of 1-quasi-hereditary algebras is not closed under Ringel-duality (ct. Theorem B).   

\begin{num}\begin{normalfont}\textbf{Theorem.} \end{normalfont} Let   $A=(KQ/\II,\ma)$   be    a 1-quasi-hereditary algebra with $1\ma i\ma n$,  then  
\begin{center}
 $R(A)$ is 1-quasi-hereditary \  if and only if  \ 
 $\displaystyle T(i)=P(1)/\left(\sum_{l\in Q_0\backslash \Lambda_{(i)}}P(l)\right)$ for any $i\in Q_0$.
\end{center} \label{R(A)}
\end{num}

Note that   the Ringel dual of a 1-quasi-hereditary algebra is 1-quasi-hereditary if and only if the equivalent conditions of Theorem~\ref{A(i)}  are satisfied. 

We now 
consider  some  properties of $R(A)$ for a 1-quasi-hereditary algebra $A$.  The vertices  in  $Q_0$ and $Q_0(R)$ will be identified. By $\ma_{(R)}$ we denote the partial order on   $Q_0(R)$, it  means $i\ma j$  if and only if $j\ma_{(R)} i$. Obviously, 
 $\left\{1\right\}=\max\left\{\left( Q_0(R),\ma_{(R)}\right)\right\}$  and  $ \left\{n\right\}=\min\left\{\left( Q_0(R),\ma_{(R)}\right)\right\}$.  For the $R(A)$-modules  we'll  use  the  index $(R)$.

\begin{num}\begin{normalfont}\textbf{Lemma.}\end{normalfont} Let $A=(KQ/\II,\ma)$ be a 1-quasi-hereditary algebra. Then,  for the  Ringel dual   $R(A)=\left(KQ_{(R)}/\II_{(R)},\ma_{(R)}\right)$ it is:
\begin{itemize}
	\item[(a)]  $P_{(R)}(n)\cong I_{(R)} (n)\cong T_{(R)}(1) $. 
	\item[(b)]   $\Delta_{(R)}(j) \hookrightarrow \Delta_{(R)}(i)$ \ if and only if  \  $\nabla_{(R)}(i) \twoheadrightarrow \nabla_{(R)}(j)$ \  if and only if  $j\ma_{(R)} i$.
\item[(c)] $\soc P_{(R)}(i) \cong \topp I_{(R)}(i) \cong  S_{(R)}(n)$ \  if and only if  \  $\soc T(i) \cong S(1)$. 
\item[(d)] $\left[\De_{(R)} (j):S_{(R)}(i)\right]= 1$    for     $i\ma_{(R)} j$,      if  \    $\topp T(i)\cong S(1)$.
\end{itemize} \label{ringeldu} \end{num}

\textit{Proof.} 
\textit{(a)}  Using   ~\ref{A0p}  and ~\ref{tilting},  we have   $ I(1)\cong T(n)$ resp.   $I_{A^{op}}(1)\cong T_{A^{op}}(n)$. By applying  $\mathscr{R}_{(A)}(-)$  resp. $\DD \left(\mathscr{R}_{(A^{op})}(-)\right)$ we obtain  $T_{(R)}(1)  \cong P_{(R)}(n) $ resp.  $T_{(R)}(1) \cong \DD(T_{R(A^{op})}(1))  \cong \DD( P_{R(A^{op})}(n))\cong I_{(R)}(n) $  
 since $R(A^{op})\cong R(A)^{op}$. Thus
  $  P_{(R)}(n)  \cong  T_{(R)}(1)  \cong  I_{(R)}(n)$.

\textit{(b)}   There exist an  exact sequence   $ \zeta : 0 \rightarrow \ger{K}\rightarrow \nabla(j) \stackrel{\pi}{\twoheadrightarrow} \nabla(i)\rightarrow 0$,  where  $\ger{K}= \Kern \pi$  if and   only if  $i\ma j$ (see ~\ref{socel}).  By applying   $\mathscr{R}_{(A)}(-)$  to    $\zeta$  we obtain an  exact sequence        $0\rightarrow \Hom_A(T,\ger{K}) \rightarrow \De_{(R)}(j) \rightarrow \De_{(R)}(i)$.
Since   $\topp \nabla(k) \stackrel{~\ref{socel}}{\cong} S(1)$ and  $[\nabla(k):S(1)]=1$ for all $k\in Q_0$ (see Sec.1 $(\ast)$),      we obtain  $\topp T \in \add S(1)$  and   $\left[\ger{K}:S(1)\right] =0$.  This implies  $\Hom_A(T,\ger{K})=0$ and  consequently  $\De_{(R)}(j) \hookrightarrow \De_{(R)}(i)$ if and only if $i\ma j$ (i.e.  $j\ma_{(R)}i$). 

The algebra $A^{op}$ is 1-quasi-hereditary, thus $\De_{(R(A^{op}))}(j)\hookrightarrow \De_{(R(A^{op}))}(i)$ if and only if $j\ma_{(R)}i$. 
Using duality,  we have  $\nabla_{(R)}(i) \twoheadrightarrow
 \nabla_{(R)}(j)$  if and only if  $j\ma_{(R)} i$.
	
\textit{(c)} \dq $\Leftarrow$\dq \   Since  $\soc T(i) \cong S(1)$, we have  $T(i)\hookrightarrow I(1)\cong T(n)$. Thus   we have an exact sequence $\xi: 0\rightarrow T(i) \rightarrow T(n)\rightarrow T(n)/T(i)\rightarrow 0$  with $T(i), T(n)$ and $ T(n)/T(i)\in \ger{F}(\nabla)$, since $\ger{F}(\nabla)$ is coresolving. Applying  $\mathscr{R}_{(A)}(-)$ to   $\xi$ yields  an exact sequence  $0\rightarrow P_{(R)}(i)\rightarrow P_{(R)}(n)\rightarrow \mathscr{R}_{(A)}(T(n)/T(i)) \rightarrow 0$. 
 Hence   \textit{(a)} implies  $\soc P_{(R)}(i)\cong S_{(R)}(n)$. 
 
According to Theorem~\ref{A(i)} and $\topp T(i)\in \add(S(1))$, we obtain that  $\soc T(i)\cong S(1) $  imlies    $\topp T(i)\cong S(1)$. Thus  $\soc \DD(T(i))\cong \soc T_{A^{op}}(i) \cong S(1)$ and therefore $\soc P_{R(A^{op})}(i) \cong S_{(R)}(n)$. Using duality, we obtain $\topp I_{(R)}(i)\cong S_{(R)}(n)$.

\dq $\Rightarrow$\dq \  Let  $\soc T(i)\cong S(1)^m$   [we know $\soc T(i)\in \add (S(1))$]. 
 Since $T(i), I(1)^m\cong T(n)^m$ and $N:=T(n)^m/T(i)\in \ger{F}(\nabla)$, applying  $ \mathscr{R}_{(A)}(-)$ to  the exact sequence $\xi: 0\rightarrow T(i)\rightarrow T(n)^m\rightarrow N \rightarrow 0$  yields  an exact  sequence 
$0\rightarrow P_{(R)}(i)\rightarrow P_{(R)}(n)^{m}\rightarrow \mathscr{R}_{(A)}(N) \rightarrow 0$.  It is sufficient to show that $P_{(R)}(n)^m \stackrel{(a)}{\cong} I_{(R)}(n)^m$  is an injective envelope of $P_{(R)}(i)$. The assumption $\soc P_{(R)}(i)\cong S_{(R)}(n)$ implies then  $m=1$ and consequently $\soc T(i)\cong S(1)$:
Assume $P_{(R)}(n)^m$ is not an injective envelope of $P_{(R)}(i)$,  then $P_{(R)}(n)$ is a direct summand of $\mathscr{R}_{(\A)}(N)$. Since $P_{(R)}(n) \stackrel{(a)}{\cong} T_{(R)}(1)$ and $(T_{(R)}(1):\De_{(R)}(1))=1$, we obtain  $\left(\mathscr{R}_{(A)}(N):\De_{(R)}(1)\right)\neq 0$.  
The properties of $\mathscr{R}_{(A)}(-)$ 
 imply   $(N:\nabla(1))=(\mathscr{R}_{(A)}(N):\De_{(R)}(1))\neq 0$.  
  Since $\left(T(n):\nabla(1)\right)\stackrel{~\ref{tilting}}{=}1$, the sequence $\xi$  provides $m=\left(T(n)^m:\nabla(1)\right)=\left(T(i):\nabla(1)\right)+\left(N:\nabla(1)\right)$. Moreover, $\left(T(i):\nabla(1)\right)\geq m$ because  $\soc T(i)\cong \nabla(1)^m$ and the filtration $0\subset \soc T(i)\subset T(i)$ can be refine to a $\nabla$-good filtration of $T(i)$. We obtain $\left(T(i):\nabla(1)\right)=m$ and therefore $\left(N:\nabla(1)\right)=0$.   We obtain a contradiction to our assumption.

\textit{(d)}  The  structure of $\De_{(R)}(j)$ yields  $[\De_{(R)}(j):S_{(R)}(i)]=\dimm_K\Hom_A(T(i), \nabla(j))$.  If   $\topp  T(i)\cong S(1)$  and  an $A$-map $F:T(i)\rightarrow \nabla(j)$ is non zero, then  $F$ is surjective    and $\dimm_K\Hom_A(T(i), \nabla(j))=1$  because $\topp \nabla(j)\cong S(1)$ and $[\nabla(j):S(1)]=1$. The
properties of $T(i)$ yield   $T(i)\twoheadrightarrow \nabla(i)$,   thus   we have a surjective map  $F':T(i)\twoheadrightarrow \nabla(i)\stackrel{~\ref{socel}}{\twoheadrightarrow} \nabla(j)$ for all $j$ with $j \ma i$.  Thus $[\De_{(R)}(j):S_{(R)}(i)]=\dimm_K\Hom_A(T(i), \nabla(j))=1$ for all $i\ma_{(R)}j$.  \hfill $\Box$
\\

\textit{Proof of the theorem.} $\dq\Rightarrow\dq$  If  $R(A)$ is 1-quasi-hereditary, then  
for any $i\in Q_0(R)$  it is  $\soc P_{(R)}(i)\cong S_{(R)}(n)$ (here $\left\{n\right\} =\min \left\{Q_0(R),\ma_{(R)}\right\}$).   Lemma~\ref{ringeldu} \textit{(d)}  implies $\soc T(i)\cong S(1)$ and  
  Theorem~\ref{A(i)}  provides $T(i)=P(1)/\left(\sum_{l\in Q_0\backslash \Lambda_{(i)}}P(l)\right)$ for any $i\in Q_0$.

$\dq\Leftarrow\dq$ If $T(i)=P(1)/\left(\sum_{j\in Q_0\backslash \Lambda_{(i)}}P(j)\right)$, then  $\soc T(i) \cong \topp T(i) \cong S(1)$   for any $i\in Q_0$ (see ~\ref{A(i)}). 
 Lemma~\ref{ringeldu} \textit{(c)} and  \textit{(b)}  provides  (3) and (4) of Definition~\ref{def1qh}.  
 
  According to Theorem~\ref{A(i)} and Lemma~\ref{factoralgebra} \textit{(a)}  the  $A$-module $T(i)$  can be considered as  the   module  $P_{(i)}(1)\cong I_{(i)}(1)$ over  a 1-quasi-hereditary algebra  $A(i)$.     Thus    $(T(i):\nabla(j))=1$ for every $j\in Q_0(i)=\Lambda_{(i)}$. We obtain $(P_{(R)}(i):\De_{(R)}(j))=1$  for every $j\in Q_0(R)$ with $i\ma_{(R)} j$ and ~\ref{ringeldu} \textit{(e)} yields  (2) of Definition~\ref{def1qh}. \hfill $\Box$
\\

If for some 1-quasi-hereditary algebra $A$ the algebra $R(A)$ is not 1-quasi-hereditary, then there exists $i\in Q_0$ such that  $\soc T(i)\cong S(1)^{m}$ with $m\geq 2$, and consequently   $(P_{(R)}(i):\De_{(R)}(1))\geq 2$.  An example of  a 1-quasi-hereditary algebra $A$ such that $R(A)$ is not  1-quasi-hereditary can be found in \cite{P}.
\\[5pt]

\textbf{Acknowledgments.}  
I would like to thank  Rolf 	Farnsteiner  and   Julian Külshammer   for  helpful   remarks and comments.

\begin{small}

\end{small}

\end{document}